\RequirePackage{ifpdf}
\ifpdf 
\documentclass[pdftex]{sigma}
\else
\documentclass{sigma}
\fi

\begin{document}
\allowdisplaybreaks

\renewcommand{\PaperNumber}{049}

\FirstPageHeading

\ShortArticleName{On One Approach to Investigation of Mechanical
Systems}

\ArticleName{On One Approach to Investigation\\ of Mechanical
Systems}

\Author{Valentin D. IRTEGOV and  Tatyana N. TITORENKO}
\AuthorNameForHeading{V.D. Irtegov and T.N. Titorenko}

\Address{Institute for Systems Dynamics and Control Theory, SB
RAS, Irkutsk, Russia}

\Email{\href{mailto:irteg@icc.ru}{irteg@icc.ru}}

\ArticleDates{Received November 18, 2005, in f\/inal form April
11, 2006; Published online May 08, 2006}

\Abstract{The paper presents some results of qualitative analysis
of Kirchhof\/f's
 dif\/ferential equations describing motion of a rigid body in ideal
 f\/luid in Sokolov's case.  The research methods are based on
 Lyapunov's  classical results. Methods of computer
 algebra implemented in the computer algebra system (CAS) ``Mathematica'' were
 also used. Combination of these methods  allowed us to obtain rather
 detailed information on qualitative properties for some classes of solutions
 of the equations.}

\Keywords{rigid body mechanics; completely integrable systems;
qualitative analysis; invariant manifolds; stability;
bifurcations; computer algebra}

\Classification{37N05; 34D20; 68W30}

\section{Introduction}
 Let us consider the problem which may be of interest both for the development of the
 method of investigation employed and for numerous applications.
 In the problem proposed for consi\-de\-ration the equations of Euler--Poisson's type
 are the model of the object under investigation.

 Euler--Poisson's dif\/ferential equations describing the motion of a
 rigid body with one f\/ixed point and their numerous
 generalizations represent one of successful mathematical models
 that is widely used in investigations of diverse physical phenomena and
 processes.

For example, the following Euler's equation for an abstract model
of an inf\/inite-dimensional dissipative top
\begin{gather}
\label{010101} \frac{d}{dt}A \psi+ \varepsilon B \psi +[ \psi, A
\psi]= \varepsilon f,
\end{gather}
where $ [ \psi, A \psi ]$ is Poisson's bracket, may be used to
describe nonsteady-state f\/lat-parallel f\/low of a viscous
incompressible f\/luid in a channel with solid walls.
 The equation describing the f\/luid motion writes:
\begin{eqnarray}
\label{020202} - \frac{\partial}{\partial t} \Delta \psi +
\varepsilon \Delta  \Delta \psi- \frac{\partial \psi}{\partial y}
\frac{\partial \Delta \psi}{\partial x} + \frac{\partial
\psi}{\partial x}  \frac{\partial \Delta \psi}{\partial y}=
\varepsilon \cos y,
\end{eqnarray}
where $ \psi (t,x,y) $ is a function of current, $ \Delta=
{\partial^{2}}/{\partial x ^{2}}
 + {\partial^{2}}/{\partial y^{2}}$ is the Laplace operator.

Under the following boundary conditions
\begin{gather*}
1) \ \ 0 < x < \frac{2 \pi }{ \alpha }, \quad   0< y < 2 \pi , \quad \alpha > 0, \\
2) \ \  \psi \left(t,x+ \frac{2 \pi}{ \alpha}, y\right)= \psi (t,x,y), \\
3) \ \ \psi |_{y=0}= \frac{\partial \psi}{\partial y}\Big|_{y=0}=
\psi\big|_{y=2 \pi}=
 \frac{\partial \psi}{\partial y}\Big|_{y=2 \pi}
\end{gather*}
(condition 3 indicates to the zero f\/luid f\/low rate) and the
following correlation between the operators
\[
A \equiv - \Delta , \qquad B \equiv  \Delta \Delta, \qquad [ \psi,
\varphi ] \equiv  \frac{\partial \psi}{\partial y}
 \frac{\partial \varphi}{\partial x}-
 \frac{\partial \psi}{\partial x} \frac{\partial \varphi}{\partial y},
\]
equation (\ref{020202}) is similar to (\ref{010101}).

Such analogies allow one to conduct, for example, analysis of
stability in problems of above type by classical methods of rigid
body dynamics, and to suggest  clear interpretation of results
obtained \cite{Oparina&2004}.

The present paper represents some results of qualitative analysis
of the dif\/ferential equations describing the motion of a rigid
body in ideal incompressible f\/luid. If the following conditions
are satisf\/ied here, i.e.\ the f\/luid possesses a single-valued
potential of rates and rests at inf\/inity, then the body motion
equations (6 ODEs) separate from the partial dif\/ferential
equations which describe the motion of f\/luid. In this case, the
motion equations of the body coincide in their form with the
corresponding Euler--Poisson equations and are called
Kirchhof\/f's equations \cite{Kir&1887,Lamb&19xx}:
\begin{gather*}
\dot M = M \times \frac{\partial H}{\partial M} + \gamma \times
\frac{\partial H}{\partial \gamma}, \qquad \dot \gamma = \gamma
\times \frac{\partial H}{\partial M},
\end{gather*}
where $M = (M_1, M_2, M_3)$, $\gamma = (\gamma_1, \gamma_2,
\gamma_3)$ are vectors of ``impulse moment'' and ``impulsive
force'', respectively.

The total kinetic energy for the body and the f\/luid writes:
\[
2T = 2 H = (A M, M) + 2 (B M, \gamma) + (C \gamma, \gamma).
\]
Here $ A$, $B$, $C$ are constant matrices. The latter are inertial
characteristics of the body and the f\/luid. By a special choice
of the origin and the direction of axes in the body it is possible
to make the matrix $ A $ diagonal, and the matrix $ C $ symmetric,
respectively.

In problems of qualitative analysis of Euler's equations, it is
possible to obtain the most complete results  when the equations
have many f\/irst integrals, for example, the equations are
completely integrable. In this case, the phase space of system of
the equations has a simple structure.

We consider Kirchhof\/f's dif\/ferential equations in Sokolov's
case \cite{Sokolov&2001}. The equations in this case represent a
completely integrable system. These write:
\begin{gather}
\dot M_1  =  M_2 M_3 + \alpha ( \gamma_2 M_1 + \gamma_1 M_2)
    + 2 \beta (\gamma_2 M_2 - \gamma_3 M_3)
    - 4 \gamma_2 \gamma_3 (2 \alpha^2 + \beta^2)
    + 4 \alpha \beta \gamma_1 \gamma_3, \nonumber \\
\dot M_2  =  4 \gamma_1 \gamma_3 (\alpha^2 + 2 \beta^2)
    - \beta (\gamma_2 M_1 + \gamma_1 M_2) -
    M_1 M_3 - 2 \alpha ( \gamma_1 M_1 - \gamma_3 M_3)
    - 4 \alpha \beta \gamma_2 \gamma_3, \nonumber \\
\dot M_3  =  4 \gamma_1 \gamma_2 (\alpha^2 - \beta^2) +
    \beta (\gamma_3 M_1 + \gamma_1 M_3) -
    \alpha (\gamma_3 M_2 + \gamma_2 M_3) -
    4 \alpha \beta (\gamma_1^2 - \gamma_2^2), \nonumber \\
\dot \gamma_1  =   \gamma_2 (2 M_3 + \alpha \gamma_1)
   + \beta (\gamma_2^2 - \gamma_3^2) - \gamma_3 M_2, \nonumber\\
\dot \gamma_2  =  - \gamma_1 (2 M_3 + \beta \gamma_2)
   - \alpha (\gamma_1^2 - \gamma_3^2) + \gamma_3 M_1, \nonumber \\
\dot \gamma_3  =  \gamma_1 (M_2 + \beta \gamma_3) -
    \gamma_2 (M_1 + \alpha \gamma_3). \label{01}
\end{gather}
Here $\alpha$, $\beta$ are arbitrary constants.

The system (\ref{01}) has the following 4 algebraic f\/irst
integrals:
\begin{gather}
 2 H  =  M_1^2 + M_2^2 + 2 M_3^2
   + 2 \alpha (\gamma_3 M_1 + \gamma_1 M_3)
   + 2 \beta (\gamma_3 M_2 + \gamma_2 M_3)
   + 4 (\beta \gamma_1 - \alpha \gamma_2)^2 \nonumber \\
\phantom{2 H  =}{} - 4 \gamma_3^2 (\alpha^2 + \beta^2) = 2 h, \nonumber \\
 V_1  =  \gamma_1 M_1 + \gamma_2 M_2 + \gamma_3 M_3 = c_1,\qquad
 V_2  =  \gamma_1^2 + \gamma_2^2 + \gamma_3^2 = c_2, \nonumber \\
 V_3  =  \big \{ 3 (\beta \gamma_1 - \alpha \gamma_2) (\beta M_1 - \alpha M_2)
 +  (2 \alpha \gamma_1 + 2 \beta \gamma_2 + M_3)
   ( (\alpha^2 + \beta^2) \gamma_3 + \alpha M_1 + \beta M_2) \big \}^2 \nonumber \\
\phantom{V_3  =}{}  + (M_3 -\alpha \gamma_1 - \beta \gamma_2)^2
   \big \{(\beta M_1 - \alpha M_2)^2 +
    (\alpha^2 + \beta^2) (2 \alpha \gamma_1
  + 2 \beta \gamma_2 + M_3)^2 \big \} = c_3. \label{02}
\end{gather}

\newpage

A lot of works were devoted to investigation of the Kirchhof\/f's
equations. A substantial part of these works is related to the
problems of integrability.

Another part of such works was devoted to the problems of
investigation of stability for permanent motions (in particular,
helical motions) of a rigid body in ideal f\/luid. The f\/irst
results in this direction go back to Lyapunov
\cite{Lyapunov&1954}. This investigation was further  developed,
and some results can be found, for example, in
\cite{Harlamov&1963,Kovalev&1968,Kozlov&2001,Ryabov&2003}.

The objective of the present work is to conduct qualitative
analysis of solutions of equations~(\ref{01}). We investigate a
class of solutions of the equations, on which the elements of
algebra of the problem's f\/irst integrals assume stationary
values. Such solutions will be called stationary (see Appendix A
for details). In particular, we  found families of stationary
solutions and families of invariant manifolds of steady motions
(IMSMs) for the system of equations under scrutiny. We  obtained
conditions of stability and instability for several families of
stationary solutions and families of IMSMs; parametric analysis of
some of these conditions was conducted. Besides, some problems of
bifurcations for both the families of stationary solutions and the
families of IMSMs branching from these solutions (in particular, a
trivial solution) were considered. Furthermore, the character of
stability for branching manifolds was taken into account.

The methods of investigation are based on classical Lyapunov's
results \cite{Lyapunov&1954,Lyapunov&1956}, in particular, on his
2nd method. Methods of computer algebra implemented in CAS
``Mathematica'' were also  used. A combination of these methods
enabled to obtain results quite interesting from our viewpoint.

\section{Obtaining stationary solutions}
Let us consider the problem of f\/inding stationary solutions and
invariant manifolds of steady motions for the system (\ref{01}).

Analysis of equations (\ref{01}) in terms of initial variables
$M_i$, $\gamma_i$, $i=1,2,3$, is rather bulky, and hence
dif\/f\/icult. Therefore, in the papers devoted to the analysis of
above equations, dif\/ferent linear transformations of variables
are applied  allowing one to reduce the equations and the
integrals to a more compact form. In the present paper we use the
following linear (not degenerate) transformation of the variables
from \cite{Borisov&2001}:
\begin{gather}
 M_1 = s_1 - \frac{1}{3} {\tilde \alpha} r_3, \qquad
 M_2 = s_2 - \frac{1}{3} {\tilde \beta} r_3, \qquad
 M_3 = s_3 + \frac{1}{3} {\tilde \alpha} r_1
  + \frac{1}{3} {\tilde \beta} r_2,\nonumber \\
 \gamma_i = r_i, \qquad \alpha = \frac{1}{3} \tilde \alpha,
 \qquad \beta = \frac{1}{3} \tilde \beta, \qquad i = 1,2,3. \label{033}
\end{gather}
The latter enabled us to f\/ind out stationary solutions and IMSMs
for the system (\ref{01}) and to perform their analysis without
going beyond standard algorithms.

 On account of the linear transformation of the variables
({\ref{033}}), the equations of motion (\ref{01}) for $\beta = 0$
will take the form:
\begin{alignat}{3}
 &  \dot r_1 = (\alpha r_1 + 2 s_3) r_2 - r_3 s_2,\qquad &&
 \dot s_1 = (\alpha r_1 + s_3) s_2 - \alpha^2 r_2 r_3,& \nonumber \\
&  \dot r_2 = r_3 s_1 - r_1 (\alpha r_1 + 2 s_3), &&
 \dot s_2 = (\alpha r_3 - s_1)(\alpha r_1 + s_3),&  \nonumber \\
&  \dot r_3 = r_1 s_2 - r_2 s_1, &&  \dot s_3 =  -\alpha r_2 s_3,
&\label{010100}
\end{alignat}
and the corresponding f\/irst integrals write:
\begin{gather}
 2 H  = (s_1^2 + s_2^2 + 2 s_3^2) + 2 \alpha r_1 s_3
   - \alpha^2 r_3^2 = 2 h, \nonumber \\
 V_1 = s_1 r_1 + s_2 r_2 + s_3 r_3 = c_1, \qquad
 V_2 = r_1^2 + r_2^2 + r_3^2 = c_2, \nonumber \\
 2 V_3 =  (\alpha r_1 s_1 + \alpha r_2 s_2 + s_1 s_3)^2
  + s_3^2 \big(s_2^2 + (\alpha r_1 + s_3)^2\big) = 2 c_3. \label{0202}
\end{gather}

We shall now consider the problem of f\/inding stationary
solutions and IMSMs for the system~(\ref{010100}). We shall apply
Routh--Lyapunov's method \cite{Rumyantsev&1995} (see also Appendix
B) for solving it. This will allow one to perform a substantial
part of computational work with the use of computer algebra
systems.

 In accordance with the method of Routh--Lyapunov, the functions
$K$ are constructed on the basis of the problem's f\/irst
integrals. We shall construct here linear combinations of these
integrals only (the combinations may be nonlinear ones):
\begin{gather}
\label{03}
 K = \lambda_0 H - \lambda_{1}V_{1}- \lambda_{2} V_{2}- \lambda_{3} V_{3},
 \qquad \lambda_i = {\rm const}.
\end{gather}
The integral $K$ represents a family of f\/irst integrals, which
is parametrized by the values $\lambda_0$, $\lambda_1$,
$\lambda_2,$ $\lambda_3$.  We enter an excessive number of
parameters into $K$ that allows us to obtain ``incomplete''
combinations of integrals by equating some part of the parameters
in $K$ to zero. It is worth to note that dif\/ferent elements of
algebra of f\/irst integrals will, generally speaking, correspond
to various stationary solutions and manifolds.

 Next, we write down stationary conditions for $K$ with respect to
all the variables $s_1$, $s_2$, $s_3$, $r_1$, $r_2$, $r_3$:
\begin{gather}
 \frac{\partial K}{\partial s_{1}} =
    \lambda_0 s_1 - \lambda_1 r_1
  - \alpha^2 \lambda_3 r_1^2 s_1
  - \alpha^2 \lambda_3 r_1 r_2 s_2
  - 2 \alpha \lambda_3 r_1 s_1 s_3
  - \alpha \lambda_3 r_2 s_2 s_3
  - \lambda_3 s_1 s_3^2 = 0,
 \nonumber \\
 \frac{\partial K}{\partial s_{2}}  =
   \lambda_0 s_2 - \lambda_1 r_2
  - \alpha^2 \lambda_3 r_1 r_2 s_1
  - \alpha^2 \lambda_3 r_2^2 s_2
  - \alpha \lambda_3 r_2 s_1 s_3
  - \lambda_3 s_2 s_3^2 = 0,
 \nonumber \\
 \frac{\partial K}{\partial s_{3}}  =
  \alpha \lambda_0 r_1 - \lambda_1 r_3
   - \alpha \lambda_3 r_1 s_1^2
   - \alpha \lambda_3 r_2 s_1 s_2
   + 2 \lambda_0 s_3
   - \alpha^2 \lambda_3 r_1^2 s_3
   - \lambda_3 s_1^2 s_3
   - \lambda_3 s_2^2 s_3
 \nonumber\\
\phantom{\frac{\partial K}{\partial s_{3}}  =}{}- 3 \alpha
\lambda_3 r_1 s_3^2
    - 2 \lambda_3 s_3^3 = 0,  \nonumber \\
 \frac{\partial K}{\partial r_{1}}  =
   \alpha \lambda_0 s_3
  - 2 \lambda_2 r_1 - \lambda_1 s_1
  - \alpha^2 \lambda_3 r_1 s_1^2
  - \alpha^2 \lambda_3 r_2 s_1 s_2
  - \alpha \lambda_3 s_1^2 s_3
  - \alpha^2 \lambda_3 r_1 s_3^2   - \alpha \lambda_3 s_3^3 = 0,\! \nonumber \\
 \frac{\partial K}{\partial r_{2}}  =
  2 \lambda_2 r_2 + \lambda_1 s_2
  + \alpha^2 \lambda_3 r_1 s_1 s_2
  + \alpha^2 \lambda_3 r_2 s_2^2
  + \alpha \lambda_3 s_1 s_2 s_3 = 0, \nonumber \\
 \frac{\partial K}{\partial r_{3}}  =
(\alpha^2 \lambda_0 + 2 \lambda_2) r_3 + \lambda_1 s_3 = 0.
\label{04}
\end{gather}
Solutions of system (\ref{04}) def\/ine stationary solutions and
IMSMs for the system (\ref{010100}) (see Appendices~A,~B). In the
general case these solutions may contain parameters $\lambda_{i}$
that appear in the family of integrals $K$ (\ref{03}), and may
hence represent a family of stationary solutions and IMSMs. Hence,
to solve the  problem formulated (obtaining stationary solutions
and IMSMs for the system (\ref{010100}) corresponding to the
family of f\/irst integrals $K$) it is necessary to solve the
system of 6 algebraic equations containing four parameters
$\lambda_0$, $\lambda_1$, $\lambda_2$, $\lambda_3$. The number of
variables is~6. The system will have a closed form if any three
f\/irst integrals (\ref{0202}) are added to it. In particular,
this means that the parameters $\lambda_i$ may be obtained as some
functions of the constants of f\/irst integrals.

All the equations of system (\ref{04}) are nonlinear, except for
the last one that allows us to slightly simplify the problem.
After removing the variable $r_3$ from the remaining equations
(\ref{04}) with the use of the last one, we obtain the following
system of 5 nonlinear algebraic equations containing 5 variables:
\begin{gather}
\lambda_0 s_1 - \lambda_1 r_1
  - \alpha^2 \lambda_3 r_1^2 s_1
  - \alpha^2 \lambda_3 r_1 r_2 s_2
  - 2 \alpha \lambda_3 r_1 s_1 s_3
  - \alpha \lambda_3 r_2 s_2 s_3
  - \lambda_3 s_1 s_3^2 = 0,
 \nonumber \\
 \lambda_0 s_2 - \lambda_1 r_2
  - \alpha^2 \lambda_3 r_1 r_2 s_1
  - \alpha^2 \lambda_3 r_2^2 s_2
  - \alpha \lambda_3 r_2 s_1 s_3
  - \lambda_3 s_2 s_3^2 = 0,
 \nonumber \\
 \alpha \lambda_0 r_1
   - \frac{\lambda_1^2}{\alpha^2 \lambda_0 + 2 \lambda_2} s_3
   - \alpha \lambda_3 r_1 s_1^2
   - \alpha \lambda_3 r_2 s_1 s_2
   + 2 \lambda_0 s_3
   - \alpha^2 \lambda_3 r_1^2 s_3
   - \lambda_3 s_1^2 s_3
 \nonumber \\
\phantom{\alpha \lambda_0 r_1}{} - \lambda_3 s_2^2 s_3
   - 3 \alpha \lambda_3 r_1 s_3^2
   - 2 \lambda_3 s_3^3 = 0,\nonumber \\
 \alpha \lambda_0 s_3
  - 2 \lambda_2 r_1 - \lambda_1 s_1
  - \alpha^2 \lambda_3 r_1 s_1^2
  - \alpha^2 \lambda_3 r_2 s_1 s_2
  - \alpha \lambda_3 s_1^2 s_3
  - \alpha^2 \lambda_3 r_1 s_3^2
  - \alpha \lambda_3 s_3^3 = 0, \nonumber \\
 2 \lambda_2 r_2 + \lambda_1 s_2
  + \alpha^2 \lambda_3 r_1 s_1 s_2
  + \alpha^2 \lambda_3 r_2 s_2^2
  + \alpha \lambda_3 s_1 s_2 s_3 = 0. \label{0404}
\end{gather}
The maximum degree of the equations belonging to the system is 3.

To the end of obtaining solutions of system (\ref{0404}), we apply
the Gr\"{o}bner bases (GB) method~\cite{Cox&2000} traditionally
used in computer algebra for solving similar systems. Software
implementation of the method can be found in many CAS. Application
of only standard tools of CAS ``Mathema\-tica'' for computing the
GB allowed us to construct the Gr\"{o}bner basis for the system
(\ref{0404}) under the following lexicographic ordering of the
variables: $r_1 > r_2 > s_2 > s_1 > s_3$. The timing for
construction of the basis, as measured on a 1100 MHz Pentium with
256 MB RAM running under Windows XP, is 2.71 seconds. Below one
can f\/ind a structure of Gr\"{o}bner's basis constructed under
the indicated ranging of the variables for the given system of
equations.
\begin{gather}
s_2 s_3  f_1 (s_1, s_3) = 0, \qquad
 \big(\big(\alpha^2 \lambda_0 + 2 \lambda_2\big) s_1
   + \alpha \lambda_1 s_3\big) f_2 (s_1, s_3) = 0, \qquad
 s_2  f_3 (s_1, s_2, s_3) = 0, \nonumber \\
 s_3  f_4 (s_1, s_3)  f_5(s_3) = 0,\qquad
 f_6 (r_2, s_1, s_2, s_3) = 0, \qquad
 s_3  f_7 (s_1, s_2, s_3) = 0, \nonumber \\
 f_8 (s_1, s_2, s_3) = 0, \qquad
 f_9 (r_1, s_1, s_3) = 0. \label{0004}
\end{gather}
Here $f_i$, $i=1,\ldots,9$, are polynomials of the variables
$s_1$, $s_2$, $s_3$, $r_1$, $r_2$. The maximum degree of the
polynomials is~7.
 These are too bulky, and so are omitted here. The system (\ref{0004})
is given in complete form in Appendix~C.

As is obvious from (\ref{0004}), the basis constructed can easily
be factorized that allows to decompose it into several subsystems,
which may be analyzed separately. Up to 12 subsystems were
identif\/ied. For each of the subsystems we  constructed a
Gr\"{o}bner basis under lexicographic ordering of the variables.
The latter enabled us to conduct some qualitative analysis of the
set of solutions of each subsystem (with respect to the
compatibility, f\/initeness or inf\/initeness of the set of the
subsystems' solutions, etc.) and hence to obtain information about
the whole set of system's (\ref{0004}) solutions (respectively,
(\ref{0404}) and (\ref{04})) and f\/ind out some groups of
solutions.

Finally, we conclude that the system (\ref{04}) has an inf\/inite
set of solutions (the variable $s_3$ is free). The following
groups of solutions (besides the trivial solution) were found out:
4~families of IMSMs and 8 families of stationary solutions.
 Some of the solutions obtained can be found below. The
solutions are given in the form representing the result of
computing.

1. The families of invariant manifolds of steady motions:
\begin{gather}
 \Bigg \{ \Bigg\{
 s_1 = -\frac{\sqrt{2 \lambda_2} (\alpha^2 \lambda_0 + 2 \lambda_2)
        \sqrt{\lambda_3 s_3^2 - \lambda_0} + \lambda_1 (2 \lambda_2
 + \alpha^2 \lambda_3 s_3^2)}{\alpha (\alpha^2 \lambda_0 + 2 \lambda_2) \lambda_3 s_3},
  \nonumber \\
   s_2 = \mp \frac{\sqrt{2 \lambda_2 z_1}}
  {\alpha (\alpha^2 \lambda_0 + 2 \lambda_2) \lambda_3 s_3},
 \quad
  r_1 = -\frac{ \sqrt{ \lambda_3 s_3^2 - \lambda_0}
    ((\alpha^2 \lambda_0 + 2 \lambda_2) \sqrt{\lambda_3 s_3^2 - \lambda_0}
    + \sqrt{2 \lambda_2} \lambda_1)}
     {\alpha (\alpha^2 \lambda_0 + 2 \lambda_2) \lambda_3 s_3}, \nonumber \\
   r_2 = \mp \frac{\sqrt{(\lambda_3 s_3^2 - \lambda_0) z_1}}
   {\alpha (\alpha^2 \lambda_0 + 2 \lambda_2) \lambda_3 s_3},
  \quad
 r_3 = -\frac{\lambda_1 s_3}{\alpha^2 \lambda_0 + 2 \lambda_2}
  \Bigg\},\nonumber\\
 \Bigg\{
 s_1 = \frac{\sqrt{2 \lambda_2} (\alpha^2 \lambda_0 + 2 \lambda_2)
        \sqrt{\lambda_3 s_3^2 - \lambda_0} - \lambda_1 (2 \lambda_2
 + \alpha^2 \lambda_3 s_3^2)}{\alpha (\alpha^2 \lambda_0 + 2 \lambda_2) \lambda_3 s_3},
   \quad
 s_2 = \mp \frac{\sqrt{2 \lambda_2 z_2}}
  {\alpha (\alpha^2 \lambda_0 + 2 \lambda_2) \lambda_3 s_3}, \nonumber \\
  r_1 = -\frac{ \sqrt{ \lambda_3 s_3^2 - \lambda_0}
    ((\alpha^2 \lambda_0 + 2 \lambda_2) \sqrt{\lambda_3 s_3^2 - \lambda_0}
    - \sqrt{2 \lambda_2} \lambda_1)}
     {\alpha (\alpha^2 \lambda_0 + 2 \lambda_2) \lambda_3 s_3}, \quad
  r_2 = \pm \frac{\sqrt{(\lambda_3 s_3^2 - \lambda_0) z_2}}
   {\alpha (\alpha^2 \lambda_0 + 2 \lambda_2) \lambda_3 s_3}, \nonumber\\
 r_3 = -\frac{\lambda_1 s_3}{\alpha^2 \lambda_0 + 2 \lambda_2}
 \Bigg\} \Bigg \}.\label{004401}
\end{gather}
 For brevity, we  introduced the following
denotations:
\begin{gather*}
 z_1 = \alpha^4 \lambda_0^3 + 4 \alpha^2 \lambda_0^2
\lambda_2
   - 2 \lambda_1^2 \lambda_2 + 4 \lambda_0 \lambda_2^2
   - (\alpha^4 \lambda_0^2 + 4 \lambda_2^2 + \alpha^2 (\lambda_1^2
   + 4 \lambda_0 \lambda_2)) \lambda_3 s_3^2\\
\phantom{z_1 =}{}
 - 2 \sqrt{2 \lambda_2} \lambda_1 (\alpha^2 \lambda_0
   + 2 \lambda_2) \sqrt{\lambda_3 s_3^2 - \lambda_0}, \\
 z_2 = \alpha^4 \lambda_0^3 + 4 \alpha^2 \lambda_0^2 \lambda_2
  - 2 \lambda_1^2 \lambda_2 + 4 \lambda_0 \lambda_2^2
  - (\alpha^4 \lambda_0^2 + 4 \lambda_2^2 + \alpha^2 (\lambda_1^2
  + 4 \lambda_0 \lambda_2)) \lambda_3 s_3^2\\
 \phantom{z_2 =}{}
  + 2 \sqrt{2 \lambda_2} \lambda_1 (\alpha^2 \lambda_0
  + 2 \lambda_2) \sqrt{\lambda_3 s_3^2 - \lambda_0}.
\end{gather*}

2. The families of stationary solutions:
\begin{gather}
 \Bigg \{ \Bigg\{
   r_1 = 0, \quad
   r_2 = \pm \frac{\sqrt{\lambda_1 \sqrt{-2 \lambda_0 \lambda_2}
      + 2 \lambda_0 \lambda_2}}{\sqrt{2} \alpha \sqrt{\lambda_2 \lambda_3}},\quad
  r_3 = 0, \quad s_1 = 0, \nonumber\\
 s_2 = \mp \frac{\sqrt{\lambda_1 \sqrt{-2 \lambda_2}
    + 2 \lambda_2 \sqrt{\lambda_0}}}
      {\alpha \sqrt{-\sqrt{\lambda_0} \lambda_3}}, \quad
 s_3 = 0 \Bigg\}, \nonumber\\
 \Bigg\{
 r_1 = \pm \frac{\sqrt{-z_3 (2 \lambda_2 + \sqrt{z_3})} (\alpha^2 \lambda_0
   + 2 \lambda_2 + \sqrt{z_3})}
 {\alpha^2 z_3 \sqrt{\lambda_3}}, \quad
 r_2 = 0,\nonumber \\
  r_3 = \pm \frac{\lambda_1 \sqrt{-z_3 (2 \lambda_2 + \sqrt{z_3})}}
  {\alpha z_3 \sqrt{\lambda_3}}, \quad
  s_1 = \pm \frac{\lambda_1
   \sqrt{-z_3 (2 \lambda_2 + \sqrt{z_3}})}{z_3 \sqrt{\lambda_3}},\quad
   s_2 = 0, \nonumber \\
 s_3 = \mp \frac{(\alpha^2 \lambda_0 + 2 \lambda_2)
       \sqrt{-z_3 (2 \lambda_2 + \sqrt{z_3})}}{\alpha z_3 \sqrt{\lambda_3}}
  \Bigg\} \Bigg \}. \label{004403}
\end{gather}
Here $
  z_3  = \alpha^4 \lambda_0^2 + 4 \lambda_2^2 +
  \alpha^2 (\lambda_1^2 + 4 \lambda_0 \lambda_2).$

Tools of computer algebra allow us to write rather easily the
obtained solutions in terms of variables of any other ``good''
coordinate system, for example, in terms of the initial variables
$M_i$, $\gamma_i$. The solutions remain stationary also in terms
of these variables.

 The families of stationary solutions (\ref{004403}) are given in terms of
the variables  $M_1$, $M_2$, $M_3$, $\gamma_1$, $\gamma_2$,
$\gamma_3$:
\begin{gather*}
 \Bigg \{ \Bigg\{
   \gamma_1 = 0, \quad
   \gamma_2 = \pm \frac{
   \sqrt{\lambda_1 \sqrt{-2 \lambda_0 \lambda_2}
      + 2 \lambda_0 \lambda_2}}{3 \sqrt{2} \alpha \sqrt{\lambda_2 \lambda_3}},
  \quad \gamma_3 = 0, \quad M_1 = 0, \\
 M_2 = \mp \frac{\sqrt{\lambda_1 \sqrt{-2 \lambda_2}
    + 2 \sqrt{\lambda_0} \lambda_2}}
      {3 \alpha \sqrt{-\sqrt{\lambda_0} \lambda_3}},  \quad
 M_3 = 0 \rbrace, \\
 \Bigg\{
 \gamma_1 = \pm \frac{\sqrt{-p_1 (2 \lambda_2 + \sqrt{p_1})} (9 \alpha^2 \lambda_0
   + 2 \lambda_2 + \sqrt{p_1})}
 {9 \alpha^2 p_1 \sqrt{\lambda_3}}, \quad
 \gamma_2 = 0,\\
   \gamma_3 = \pm \frac{\lambda_1 \sqrt{-p_1 (2 \lambda_2 + \sqrt{p_1})}}
  {3 \alpha p_1 \sqrt{\lambda_3}},\quad
 M_1 = \pm \frac{2 \lambda_1
   \sqrt{-p_1 (2 \lambda_2 + \sqrt{p_1})}}{3 p_1 \sqrt{\lambda_3}},\quad
   M_2 = 0, \\
 M_3 = \mp \frac{(18 \alpha^2 \lambda_0 + 4 \lambda_2 - \sqrt{p_1})
       \sqrt{-p_1 (2 \lambda_2 + \sqrt{p_1})}}{9 \alpha p_1 \sqrt{\lambda_3}}
  \Bigg\} \Bigg \}.
\end{gather*}
Here $  p_1 = 81 \alpha^4 \lambda_0^2 + 4 \lambda_2^2 +
 9 \alpha^2 (\lambda_1^2 + 4 \lambda_0 \lambda_2)$.

 Analysis of IMSMs (\ref{004401}) showed that, after transforming these
expressions to the form, which does not contain the problem's
variables under radicals, we obtain one family of IMSMs. The
system of equalities (\ref{004401}) is a representation of the
latter in dif\/ferent maps only. The equations, which def\/ine the
family of IMSMs, write:
\begin{gather}
 \alpha^4 \lambda_1^2 \lambda_3^2 s_3^4
   + 2 a_{11} \alpha^3 \lambda_1 \lambda_3^2 s_1 s_3^3
   + a_{11}^2 \alpha^2 \lambda_3^2 s_1^2 s_3^2
   - 2 a_{10} \lambda_2 \lambda_3 s_3^2
   + 4 a_{11} \alpha \lambda_1 \lambda_2 \lambda_3 s_1 s_3
  + 2 a_6 \lambda_2 = 0, \nonumber \\
 a_{11}^4 \alpha^4 \lambda_3^4 s_2^2 s_3^4
   + 4 a_{11}^2 a_9 \alpha^2 \lambda_2 \lambda_3^3 s_2 s_3^4
   + 4 a_9^2 \lambda_2^2 \lambda_3^2 s_3^4
   - 4 a_{11}^2 a_8 \alpha^2 \lambda_2 \lambda_3^2 s_2 s_3^2\nonumber\\
 \qquad{} - 8 a_4 \lambda_2^2 \lambda_3 s_3^2
  + 4 a_6^2 \lambda_2^2 = 0, \nonumber \\
 a_{11}^2 \lambda_3^2 s_3^4
   + 2 a_{11}^2 \alpha \lambda_3^2 r_1 s_3^3
   + a_{11}^2 \alpha^2 \lambda_3^2 r_1^2 s_3^2
   - 2 a_7 \lambda_3 s_3^2
   - 2 a_{11}^2 \alpha \lambda_0 \lambda_3 r_1 s_3
   + a_6 \lambda_0 = 0, \nonumber\\
 a_9^2 \lambda_3^4 s_3^8
   + 2 a_{11}^2 a_9 \alpha^2 \lambda_3^4 r_2 s_3^6
   - 2 a_3 \lambda_3^3 s_3^6
   + a_{11}^4 \alpha^4 \lambda_3^4 r_2^2 s_3^4
   - 2 a_{11}^2 a_5 \alpha^2 \lambda_3^3 r_2 s_3^4
   + a_1 \lambda_3^2 s_3^4 \nonumber \\
\qquad{}  + 2 a_{11}^2 a_8 \alpha^2 \lambda_0 \lambda_3^2 r_2
s_3^2
   - 2 a_2 \lambda_0 \lambda_3 s_3^2
   + a_6^2 \lambda_0^2 = 0, \nonumber \\
 a_{11} r_3 + \lambda_1 s_3 = 0, \label{004402}
\end{gather}
where $a_i$, $i=1,\ldots,11$, are polynomials of the parameters
$\lambda_0$, $\lambda_1$, $\lambda_2$, $\lambda_3$. These are too
cum\-ber\-some, and so are omitted herein.

From the geometric viewpoint, equations (\ref{004402}) -- for each
f\/ixed collection of parameters of the family $\lambda_0$,
$\lambda_1$, $\lambda_2$, $\lambda_3$ -- describe the curves lying
in ${\mathbb R}^6$ at the intersection of the three 4th-order
hypersurfaces, the 8th-order hypersurface and the hyperplane.

The complete analysis of the family of IMSMs (\ref{004402}) and of
motions on it is not given in the present paper. This analysis is
rather nontrivial and may be a subject matter for another paper.

\section{Investigation of stability for stationary solutions\\ and for
IMSMs}

Let us consider the problem of stability for a series of families
of stationary solutions and families of IMSMs of system
(\ref{010100}).

\subsection{Investigation of invariant manifolds}
Consider equations (\ref{04}) under the condition $ \lambda_{0}=
\lambda_{1}= \lambda_{3}=0 $. It can be readily seen that thy have
the following solution $ r_{1}= r_{2}=r_{3}=0 $, which def\/ines
the 3-dimensional invariant manifold of steady motions for the
equations (\ref{010100}). The latter allow one to def\/ine the
vector f\/ield on the IMSM (i.e.\ to reduce the initial system of
equations to the IMSM):
\begin{gather}
\label{11331} \dot {s_{1}}= s_{2} s_{3}, \qquad \dot {s_{2}}=-
s_{1} s_{3}, \qquad \dot{s_{3}}= 0.
\end{gather}
If we consider the initial values of the variables $ s_{i}$,
$i=1,2,3$, as parameters, it is possible to assume that here we
deal with the family of invariant manifolds. Hence, for $
s_{3}^{0}=0 $ this family adjuncts to the zero solution of the
problem, i.e.\ it has at least one common point with the zero
solution.

 When using the f\/irst integral $ V_{2}= r_{1}^{2}+ r_{2}^{2}+
r_{3}^{2}$, which assumes the minimum value on the family of
IMSMs, we easily conclude on stability of the elements of the
family of IMSMs~(\ref{11331}).

The vector f\/ield on the elements of the family of IMSMs has the
two f\/irst integrals:
\[
 W_{1}= s_{1}^{2}+ s_{2}^{2}=m_{1}, \qquad  W_{2}=s_{3}=m_{2}.
\]
 Consequently, in this case we may speak of investigation of the 2nd-level
 stationary solutions (see Appendix~A). Finding such solutions and their
 analysis are trivial in the computational aspect, and we will not concentrate
 on them here. Let us consider a more complex case.

\subsection{Investigation of 2nd-level stationary solutions}

 By constructing the Gr\"{o}bner bases with respect to the problem's variables
and to one or several parameters $\lambda_0$, $\lambda_1$,
$\lambda_2$, $\lambda_3$ for the subsystems of system (\ref{0004})
we can also obtain solutions under some conditions imposed on
above parameters. For example, we  constructed the Gr\"{o}bner
basis with respect to the variables $r_1$, $r_2$, $s_1$ and the
parameter $\lambda_2$ for the subsystem
\begin{gather}
s_2 s_3 = 0, \qquad
 \big((\alpha^2 \lambda_0 + 2 \lambda_2) s_1
   + \alpha \lambda_1 s_3\big) \ f_2 (s_1, s_3) = 0, \qquad
s_2 = 0, \qquad
s_3 = 0, \nonumber \\
 f_6 (r_2, s_1, s_2, s_3) = 0, \qquad
 f_8 (s_1, s_2, s_3) = 0, \qquad
 f_9 (r_1, s_1, s_3) = 0 \label{005511}
\end{gather}
of system (\ref{0004}). It enabled us to obtain solutions of the
system under the following conditions imposed on the parameter
$\lambda_2$:
\begin{gather}
\label{556611}
 \lambda_2 = -\frac{2 \alpha^2 \lambda_0^2 +
\lambda_1^2
 + \lambda_1 \sqrt{4 \alpha^2 \lambda_0^2 + \lambda_1^2}}{4
 \lambda_0},
\\
\label{556622}  \lambda_2 = -\frac{2 \alpha^2 \lambda_0^2 +
\lambda_1^2
 - \lambda_1 \sqrt{4 \alpha^2 \lambda_0^2 + \lambda_1^2}}{4 \lambda_0}.
\end{gather}
The solutions obtained can be found in Appendix~D. The solutions
are given right in the form representing the result of computing.

Likewise in case of (\ref{004401}), after transforming the
expressions of solutions obtained to the form, which does not
contain the problem's variables under radicals, we  found out that
these solutions represent the families of IMSMs for the system
(\ref{010100}), these solutions being written in terms of the maps
for these families. Finally, we have 2 families of IMSMs.

The f\/irst family of IMSMs, which corresponds to $\lambda_2$
(\ref{556611}), can be written as:
\begin{gather}
 2 \alpha^4 \lambda_0^2 r_1^2
 + \left(2 \alpha^2 \lambda_0^2 + \lambda_1
\left(\lambda_1 - \sqrt{4 \alpha^2 \lambda_0^2 + \lambda_1^2}\right)\right) s_2^2 \nonumber\\
\qquad{} =
 \left(2 \alpha^2 \lambda_0^2 + \lambda_1 \left(\lambda_1
 - \sqrt{4 \alpha^2 \lambda_0^2 + \lambda_1^2}\right)\right) \frac{\lambda_0}{\lambda_3},
 \nonumber \\
 s_1^2 + s_2^2 = \frac{\lambda_0}{\lambda_3}, \qquad
 s_3 = 0, \qquad r_3 = 0, \qquad
 2 \alpha^2 \lambda_0 r_2 + \left(\lambda_1 - \sqrt{4 \alpha^2 \lambda_0^2
 + \lambda_1^2}\right) s_2 = 0,\label{1133101}
\end{gather}
and the second family of IMSMs, which corresponds to $\lambda_2$
(\ref{556622}), can be written as:
\begin{gather}
 2 \alpha^4 \lambda_0^2 r_1^2
 + \left(2 \alpha^2 \lambda_0^2 + \lambda_1
\left(\lambda_1 + \sqrt{4 \alpha^2 \lambda_0^2 + \lambda_1^2}\right)\right) s_2^2\nonumber\\
\qquad{} =
 \left(2 \alpha^2 \lambda_0^2 + \lambda_1 \left(\lambda_1
 + \sqrt{4 \alpha^2 \lambda_0^2 + \lambda_1^2}\right)\right) \frac{\lambda_0}{\lambda_3},
 \nonumber \\
 s_1^2 + s_2^2 = \frac{\lambda_0}{\lambda_3}, \qquad
 s_3 = 0, \qquad r_3 = 0, \qquad
 2 \alpha^2 \lambda_0 r_2 + \left(\lambda_1 + \sqrt{4 \alpha^2 \lambda_0^2
 + \lambda_1^2}\right) s_2 = 0.\label{1133102}
\end{gather}
From the geometric viewpoint, each of these families of IMSMs --
for each f\/ixed collection of parameters of the family
$\lambda_0$, $\lambda_1$, $\lambda_3$ -- describes the curves
lying in ${\mathbb R}^6$ at the intersection of the three
hyperplanes, one elliptic ``cylinder'' and one circular
``cylinder''.

The vector f\/ield on elements of the family of IMSMs
(\ref{1133101})
 is given by the dif\/ferential equation:
\begin{gather}
\label{11330}
 \dot s_2 = -\frac{2 \alpha \lambda_0 (\lambda_0 - \lambda_3 s_2^2)}
     {\big(\lambda_1 + \sqrt{4 \alpha^2 \lambda_0^2 + \lambda_1^2}\big) \lambda_3},
\end{gather}
which is derived from equations (\ref{010100}) after removing
$s_1$, $s_3$, $r_1$, $r_2$, $r_3$ from them with use of
expressions (\ref{1133101}).

The vector f\/ield on elements of the family of IMSMs
(\ref{1133102}) is given by the dif\/ferential equation:
\begin{gather}
\label{11440}
 \dot s_2 = -\frac{2 \alpha \lambda_0 (\lambda_0 - \lambda_3 s_2^2)}
     {(\lambda_1 - \sqrt{4 \alpha^2 \lambda_0^2 + \lambda_1^2}) \lambda_3},
\end{gather}
which can also be derived from equations (\ref{010100}) after
removing $s_1$, $s_3$, $r_1$, $r_2$, $r_3$ from them with the use
of~(\ref{1133102}).

Consider the problem of f\/inding 2nd-level stationary solutions
for the system (\ref{11330}) and investigation of their stability
on elements of the family of IMSMs (\ref{1133101}).

As obvious from (\ref{11330}), solutions of the form
\begin{gather}
\label{113300}
 s_2^0 = -\frac{\sqrt{\lambda_0}}{\sqrt{\lambda_3}}, \qquad
 s_2^0 = \frac{\sqrt{\lambda_0}}{\sqrt{\lambda_3}}
\end{gather}
are its stationary solutions. Now we investigate their stability
by Lyapunov's method \cite{Lyapunov&1956}. Let us consider the
f\/irst of these solutions (\ref{113300}).

According to the above method, we consider the solution $s_2^0 =
-\sqrt{\lambda_0}/\sqrt{\lambda_3}$ in the capacity of undisturbed
one. Next, we introduce the deviations $z = s_2 - s_2^0$ of
disturbed motion from undisturbed one, and consider the function $
 V = \frac{1}{2} z^2
$ in the capacity of the Lyapunov function.

The equation of disturbed motion writes:
\begin{gather}
\label{111111} \dot z =
  - \frac{2 \alpha \lambda_0 z^2 }
   {\lambda_1 + \sqrt{4 \alpha^2 \lambda_0^2  + \lambda_1^2}}
  - \frac{4 \alpha \lambda_0^{3/2} z }
   {\big(\lambda_1 + \sqrt{4 \alpha^2 \lambda_0^2  + \lambda_1^2}\big) \lambda_3}.
\end{gather}
The derivative of the function $V$ due to the dif\/ferential
equation (\ref{111111})  up to the 2nd order terms writes:
\begin{gather*}
 \dot V = - \frac{4 \alpha \lambda_0^{3/2} z^2 }
   {\big(\lambda_1 + \sqrt{4 \alpha^2 \lambda_0^2  + \lambda_1^2}\big) \lambda_3} +
   \bar V_n,
\end{gather*}
where $\bar V_n$ is the $n$-order terms ($n > 2$).

The function $V$ is positive-def\/inite. From the expression for
its derivative it is obvious that it will be negative-def\/inite
when the following conditions imposed on the parameters
$\lambda_i$ hold:
\begin{gather}
\label{005500}
 \alpha > 0 \vee \lambda_0 > 0 \vee (\lambda_1 < 0 \vee \lambda_3 > 0
   \wedge \lambda_1 > 0 \vee \lambda_3 > 0).
\end{gather}
According to the Lyapunov theorem \cite{Lyapunov&1956} on
stability of undisturbed motion, satisfaction of these conditions
means that the solution investigated is asymptotically stable on
elements of the family of IMSMs (\ref{1133101}).

According to the above theorem, the solution under scrutiny is
unstable on elements of the family of IMSMs (\ref{1133101})  when
the following conditions hold:
\begin{eqnarray}
\label{0055000}
 \alpha < 0 \vee \lambda_0 > 0 \vee (\lambda_1 < 0 \vee \lambda_3 > 0
   \wedge \lambda_1 > 0 \vee \lambda_3 > 0).
\end{eqnarray}

Substitution of the stationary solution $s_2 = -
\sqrt{\lambda_0}/\sqrt{\lambda_3}$ into (\ref{1133101}) allows one
to obtain the solution corresponding to it in the whole space of
variables $s_1$, $s_2$, $s_3$, $r_1$, $r_2$, $r_3$:
\begin{gather}
\label{222222}
 \left\{ r_1 = 0, \
  r_2 = \frac{\lambda_1 - \sqrt{4 \alpha^2 \lambda_0^2 + \lambda_1^2}}
    {2 \alpha^2 \sqrt{\lambda_0 \lambda_3}}, \
    r_3 = 0, \ s_1 = 0, \
   s_2 = - \frac{\sqrt{\lambda_0}}{\sqrt{\lambda_3}}, \
   s_3 = 0
  \right\}.
\end{gather}
Analysis of (\ref{222222}) showed that this solution is unstable
in the sense of Lyapunov \cite{Lyapunov&1954} under the following
conditions imposed on the parameters $\lambda_i$:
\begin{gather}
 \alpha < 0 \vee (\lambda_0 < 0 \vee \lambda_1 > 0 \vee \lambda_3 < 0
   \wedge \lambda_0 > 0 \vee \lambda_1 > 0 \vee \lambda_3 > 0)\nonumber \\
 \wedge \
 \alpha > 0 \vee (\lambda_0 < 0 \vee \lambda_1 > 0 \vee \lambda_3 < 0
   \wedge \lambda_0 > 0 \vee \lambda_1 > 0 \vee \lambda_3 > 0). \label{005501}
\end{gather}
There are roots having  positive real part among the roots of the
characteristic equation constructed for the equations
(\ref{010100}) linearized in the neighbourhood of the solution
(\ref{222222}) when conditions (\ref{005501}) hold.

Comparison of the conditions (\ref{005500}) and (\ref{005501})
shows that under the same conditions imposed on the parameters
$\lambda_i$:
\[
 \alpha > 0 \vee \lambda_0 > 0 \vee \lambda_1 > 0 \vee \lambda_3 > 0,
\]
the f\/irst of the solutions (\ref{113300}) asymptotically stable
on the IMSM (\ref{1133101}) corresponds to the solution which is
unstable in the whole space of the problem's variables.

Similarly, comparison of the conditions (\ref{0055000}) and
(\ref{005501}) allows us to conclude: when the following
conditions
\[
 \alpha < 0 \vee \lambda_0 > 0 \vee \lambda_1 > 0 \vee \lambda_3 > 0
\]
hold, the solution unstable on the IMSMs (\ref{1133101})
corresponds to the solution unstable in the whole space of the
problem's variables.

We also investigated stability of the 2nd solution of
(\ref{113300}) and of stationary solutions of
system~(\ref{11440}). The results appeared to be similar to those
given above.

As far as behaviour of the IMSMs themselves (\ref{1133101}),
(\ref{1133102}) is concerned, note the following. The manifolds
intersect and coincide under the condition $\lambda_0 = \lambda_1
= 0$. As a result, we obtain the invariant manifold def\/ined by
the equations:
\begin{gather}
\label{221100}
 s_1 = 0, \qquad s_2 = 0, \qquad s_3 = 0, \qquad r_3 = 0.
\end{gather}
The vector f\/ield on IMSM (\ref{221100}) is given by the
dif\/ferential equations:
\[
\dot r_1 = \alpha r_1 r_2, \qquad
 \dot r_2 = -\alpha^2 r_1^2.
\]

\section{On stability of the zero solution and of IMSMs adjunct to it}
   Let us consider a class of special stationary solutions. These solutions
possess the following properties: several of the problem's f\/irst
integrals assume stationary values on the solutions and, as
a~rule, there is a bifurcation of stationary solutions (invariant
manifolds) of various dimensions in their neighbourhood.

 For example, in completely integrable cases of the problem of rigid
body's motion having one f\/ixed point, when all the f\/irst
integrals are quadratic, the manifolds of dimension~3 (it is half
of the number of variables used in describing the problem) are
typical IMSMs branching from special permanent rotations.
Furthermore, in many cases there is the following relationship
between the stability of special permanent rotations and the
stability of IMSMs branching from them: {\it branching of stable
IMSMs from special permanent rotations represents the necessary
and sufficient stability condition for them}. If algebraic
f\/irst integrals of the problem under consideration are not only
quadratic then the relationship between the property of stability
of special stationary solutions and the property of stability of
IMSMs branching from them is more complex.

As far as Kirchhof\/f's equations are concerned that is obvious
from their analysis in Sokolov's case, the IMSMs of both even and
odd dimensions can branch from special stationary solutions
(helical motions). Now, we consider a particular example of such a
bifurcation of stationary solutions.

Consider the problem of stability of both the zero solution and
the IMSMs adjunct to it, i.e.\ the stability of invariant
manifolds which have at least one common point with the zero
solution. Equality of the system's (\ref{04}) Jacobian, which is
computed for zero values of the variables $s_i$, $r_i$
$i=1,\ldots,3$, to zero is the condition of existence of such
IMSMs.

 The Jacobian of (\ref{04}) for
 $s_1 = 0$, $s_2 = 0$, $s_3 = 0$, $r_1 = 0$, $r_2 = 0$, $r_3 = 0$ writes:
\[
 J = \big(\lambda_1^2 + 2 \lambda_0 \lambda_2\big)
 \big(\alpha^2 \lambda_0^2 + \lambda_1^2 + 2 \lambda_0 \lambda_2\big)
 \big(\alpha^2 \lambda_0^2 + \lambda_1^2 + 4 \lambda_0 \lambda_2\big).
\]
The following possibilities for satisfaction of the equality $J =
0$ were considered (when the zero solution is special):
\begin{gather}
1) \ \lambda_0 = 0, \ \lambda_1 = 0; \qquad
 2) \ \lambda_1 = 0, \ \lambda_2 = 0; \qquad
 3) \ \lambda_2 = -\frac{\lambda_1^2}{2 \lambda_0};  \nonumber\\
 4) \ \lambda_2 = -\frac{\alpha^2 \lambda_0^2 + \lambda_1^2}{2 \lambda_0}; \qquad
 5) \ \lambda_2 = -\frac{\alpha^2 \lambda_0^2 + \lambda_1^2}{4 \lambda_0}.\label{5500}
\end{gather}
Obtaining solutions of stationary equations (\ref{04}) for the
indicated values of $\lambda_2$ (\ref{5500}) with the aid of
Gr\"{o}bner's bases technique  allowed us to f\/ind out a
suf\/f\/iciently large number of families of IMSMs adjunct to the
zero solution. Some of these solutions are adduced below:
\begin{gather}
\label{117}
 \left\{
 s_1 = -\frac{\lambda_1}{\alpha \lambda_0} s_3, \ s_2 = 0,\
 r_1 = -\frac{s_3}{\alpha}, \ r_2 = 0, \ r_3 = \frac{\lambda_0}{\lambda_1} s_3
  \right\} \qquad {\rm for}\quad
 \lambda_2 = -\frac{\alpha^2 \lambda_0^2 + \lambda_1^2}{2 \lambda_0};\\
 \left\{
 s_3 = 0, \ r_1 = 0, \ r_2 = 0, \ r_3 = 0
  \right\} \qquad {\rm for} \quad
 \lambda_0 = 0,  \ \ \lambda_1 = 0; \nonumber \\
 \left\{
 s_1 = 0, \ s_2 = 0, \ s_3 = 0, \ r_1 = 0, \ r_3 = 0
  \right\} \qquad {\rm for} \quad
 \lambda_1 = 0, \ \ \lambda_2 = 0; \nonumber \\
 \left\{
 s_1 = 0, \ s_2 = \frac{\lambda_1}{\lambda_0} r_2, \ s_3 = 0, \
 r_1 = 0, \ r_2 = 0
  \right\} \qquad {\rm for} \quad
 \lambda_2 = -\frac{\lambda_1^2}{2 \lambda_0}, \ \  \lambda_3 = 0; \nonumber \\
 \left\{
 s_1 = -\frac{2 \alpha \lambda_0 \lambda_1}
  {\alpha^2 \lambda_0^2 - \lambda_1^2} s_3, \ s_2 = 0, \
 r_1 = -\frac{2 \alpha \lambda_0^2}{\alpha^2 \lambda_0^2 - \lambda_1^2} s_3, \
 r_2 = 0, \  r_3 = -\frac{2 \alpha \lambda_0 \lambda_1}
  {\alpha^2 \lambda_0^2 - \lambda_1^2} s_3
  \right\} \nonumber\\
 \qquad {\rm for} \quad
\lambda_2 = -\frac{\alpha^2 \lambda_0^2 + \lambda_1^2}{4
\lambda_0}, \ \
\lambda_3 = 0; \label{1172} \\
\left\{
 s_1 = 0, \
 s_2 = 0, \
 s_3 = -\frac{\alpha} {2} r_1, \
 r_2 = 0, \ r_3 = 0\right\} \nonumber\\
 \qquad {\rm for}\quad
\lambda_2 = -\frac{\alpha^2 \lambda_0^2 + \lambda_1^2}{4
\lambda_0}, \ \ \lambda_1 = 0, \  \ \lambda_3 = 0. \nonumber
\end{gather}
These solutions represent the simplest of the invariant manifolds
obtained: the zero vector f\/ield is def\/ined on all the IMSMs
(\ref{117}), (\ref{1172}). The latter attributes some
specif\/icity to the problem of stability analysis for such
manifolds.

We use the method of Lyapunov functions, in particular, the
Routh--Lyapunov method, to investigate stability of the zero
solution as well as the families of IMSMs (\ref{117}),
(\ref{1172}) adjunct to~it.

\subsection{Investigation of stability of the zero solution}

The procedure of obtaining suf\/f\/icient conditions of stability
for the stationary solutions by Routh--Lyapunov's method is
practically reduced to the verif\/ication (in the simplest case)
of signdef\/initeness of the 2nd variation of the integral $K$
(\ref{03}) in the neighbourhood of the statio\-nary solution under
scrutiny. The suf\/f\/icient conditions can be made ``softer'', if
signdef\/initeness of $\delta^2 K$ is considered on the manifold
def\/ined by the f\/irst variations of each of $ m-1$ integrals
(where $m$ is the number of vanishing integrals in  $K$).

Obtaining suf\/f\/icient stability conditions for the zero
solution by this technique is rather trivial in the computational
aspect. Note only that it is stable in the sense of Lyapunov when
the following restrictions imposed on the problem's parameters are
satisf\/ied:
\begin{gather}
 \alpha > 0 \vee \lambda_0 > 0 \vee
\left( \lambda_1 > 0 \vee
    \lambda_2 < -\frac{\alpha^2 \lambda_0^2 + \lambda_1^2}{2 \lambda_0} \wedge
   \lambda_1 < 0 \vee
   \lambda_2 < -\frac{\alpha^2 \lambda_0^2 + \lambda_1^2}{2 \lambda_0}\right ) \wedge
   \nonumber \\
 \alpha < 0 \vee \lambda_0 > 0 \vee
 \left ( \lambda_1 > 0 \vee
   \lambda_2 < -\frac{\alpha^2 \lambda_0^2 + \lambda_1^2}{2 \lambda_0} \wedge
  \lambda_1 < 0 \vee
   \lambda_2 < -\frac{\alpha^2 \lambda_0^2 + \lambda_1^2}{2 \lambda_0} \right ).
   \label{060606}
\end{gather}
It is possible to slightly ``soften'' the stability conditions
obtained, considering that the zero solution is special. To this
end, we solve the problem of choosing the ``best'' f\/irst
integral, i.e.\ the one which gives the ``most soft''
suf\/f\/icient stability conditions. The solution of this problem
can be obtained, for example, with the use of parametric analysis
of stability conditions, i.e.\ their minimization with respect to
one or several parameters.

Consider one of the conditions imposed on the parameters
(\ref{060606}):
\begin{gather}
\label{070707} \lambda_2 < -\frac{\alpha^2 \lambda_0^2 +
\lambda_1^2}{2 \lambda_0} \qquad {\rm or} \qquad \Lambda = 2
\lambda_2 \lambda_0 + \alpha^2 \lambda_0^2 - \lambda_1^2 < 0.
\end{gather}
Next, we f\/ind a stationary value $\Lambda$ with respect to
$\lambda_0$:
\[
 \frac{\partial \Lambda}{\partial \lambda_0} =
 2 \lambda_2 + 2 \alpha^2 \lambda_0 = 0.
\]
The latter holds for $\lambda_0 = -\lambda_2/\alpha^2$. Its
substitution into (\ref{070707}) gives the following condition
imposed on the parameters:
\[
 -\frac{\lambda_2^2}{\alpha^2} - \lambda_1^2 < 0.
\]
Hence we have stability of the equilibrium position under $\alpha
\neq 0$ and $\lambda_1 \neq 0$ or $\alpha \neq 0$ and $\lambda_2
\neq 0$.

\subsection{Investigation of stability of IMSMs adjunct to the zero solution}

Let us investigate stability of the family of IMSMs (\ref{117}).
The vector f\/ield on the elements of the family of IMSMs
(\ref{117}) is given by the dif\/ferential equation
\begin{gather}
\label{118} \dot s_3 = 0,
\end{gather}
derived from equations (\ref{010100}) after removing $s_1$, $s_2$,
$r_1$, $r_2$, $r_3$ from them with the aid of
expressions~(\ref{117}).

Using (\ref{118}), we may conclude that the elements of the family
of IMSMs obtained represent some curves in ${\mathbb R}^6$, over
each point of which the one-dimensional family of solutions ($s_3
= s_3^0 = {\rm const}$) for the equation (\ref{118}) is def\/ined.
Such IMSMs with the ``bundle'' def\/ined on them will be called
``framed invariant manifolds''. Each point in the framed IMSM
corresponds to some helical motion of a rigid body.

To the end of obtaining suf\/f\/icient stability conditions for
framed IMSMs we use the standard Lyapunov's technique.

The second variation of $K$ in the neighbourhood of some helical
motion $s_3^0$, which lies on the chosen IMSM represented in terms
of deviations
\begin{gather*}
 z_1 = r_{1}+ \frac {s_{3}}{\alpha}, \qquad z_2 = r_{2}, \qquad z_3 = r_{3}- \frac
{\lambda_{0}s_{3}}{\lambda_{1}}, \qquad z_4 =s_{1}+
\frac{\lambda_{1}s_{3}}
 { \alpha \lambda_{0}}, \\ z_5 = s_{2}, \qquad z_6 = s_3 - s_3^0
\end{gather*}
writes:
\begin{gather*}
 \delta^2 K  =  \frac{\big(\alpha^2 \lambda_0^2 + \lambda_1^2\big) \big( \lambda_0
  - \lambda_3 {s_3^0}^2\big)}{2 \lambda_0^2} z_1^2
  + \frac{\alpha^2 \lambda_0^2 + \lambda_1^2}{2 \lambda_0} z_2^2
  + \frac{\lambda_1^2}{2 \lambda_0} z_3^2
  - \lambda_1 z_1 z_4
  + \frac{\lambda_0}{2} z_4^2 \\
\phantom{\delta^2 K  =}{}  - \lambda_1 z_2 z_5
   + \frac{1}{2}\big(\lambda_0 - \lambda_3 {s_3^0}^2\big) z_5^2.
\end{gather*}
And the respective variations of the f\/irst integrals $H, V_1,
V_2$ are:
\begin{gather}
 \delta H = \left(\alpha z_1
     - \frac{\alpha^2 \lambda_0}{\lambda_1} z_3
     - \frac{\lambda_1}{\alpha \lambda_0} z_4
     + \frac{\lambda_1^4
     - \alpha^4 \lambda_0^4}{\alpha^2 \lambda_0^2 \lambda_1^2} z_6\right) s_3^0 = 0,
   \nonumber \\
 \delta V_1 = \left(-\frac{\lambda_1}{\alpha \lambda_0} z_1
     + z_3 - \frac{1}{\alpha} z_4
     + \frac{2 (\alpha^2 \lambda_0^2 + \lambda_1^2)}
       {\alpha^2 \lambda_0 \lambda_1} z_6\right) s_3^0 = 0, \nonumber\\
 \delta V_2 = 2 \left(-\frac{1}{\alpha} z_1
     + \frac{\lambda_0}{\lambda_1} z_3
     + \frac{\alpha^2 \lambda_0^2 + \lambda_1^2}{\alpha^2 \lambda_1^2} z_6\right) s_3^0 = 0.\label{776655}
\end{gather}
Conditions of signdef\/initeness for $\delta^2 K$ are
suf\/f\/icient stability conditions for the elements of the family
of IMSMs (\ref{117}). After trivial transformations, these
conditions may be written as follows:
\begin{gather*}
 \lambda_0 - \lambda_3 {s_3^0}^2 > 0, \qquad
 \alpha^2 \lambda_0^3 - \big(\alpha^2 \lambda_0^2
  + \lambda_1^2\big) \lambda_3 {s_3^0}^2 > 0, \\
 \lambda_1^3 \lambda_3 {s_3^0}^2
  + \alpha^2 \lambda_0^2 \lambda_1 \big(\lambda_3 {s_3^0}^2 - \lambda_0\big) \neq 0, \qquad
 \lambda_0 > 0, \qquad \lambda_1 \neq 0.
\end{gather*}

Now we extend the problem and investigate  stability of a body's
helical motions corresponding to the elements of the framed IMSM.

Assuming that $s_3^0 \neq 0$, we eliminate the variables $z_1$,
$z_4$ with the use of equations (\ref{776655}) (because there are
only two linearly independent ones) from $\delta^2 K$. As a
result, the following quadratic form yields:
\begin{gather*}
 \delta^2 \tilde K  =
  \frac{\alpha^2 \lambda_0^2 + \lambda_1^2}{2 \lambda_0} z_2^2
  + \frac{\lambda_1^4 + \alpha^2 \lambda_0 \big(\alpha^2 \lambda_0^2
  + \lambda_1^2\big)\big(\lambda_0 - \lambda_3 {s_3^0}^2\big)}{2 \lambda_0 \lambda_1^2} z_3^2
 - \lambda_1 z_2 z_5
 + \frac{1}{2} \big(\lambda_0 - \lambda_3 {s_3^0}^2\big) z_5^2 \\
\phantom{\delta^2 \tilde K =}{} + \frac{(\alpha^2 \lambda_0^2 +
\lambda_1^2) \big (\alpha^2 \lambda_0^3 - (\alpha^2 \lambda_0^2 +
\lambda_1^2)
 \lambda_3 {s_3^0}^2\big)}{\lambda_0 \lambda_1^3} z_3 z_6\\
 \phantom{\delta^2 \tilde K =}{}
 + \frac{(\alpha^2 \lambda_0^2 + \lambda_1^2)^2
 \big(\alpha^2 \lambda_0^3 - (\alpha^2 \lambda_0^2 + \lambda_1^2)
 \lambda_3 {s_3^0}^2\big)}{2 \alpha^2 \lambda_0^2 \lambda_1^4} z_6^2.
\end{gather*}

Conditions of signdef\/initeness for $\delta^2 \tilde {K}$ are
suf\/f\/icient stability ones for the helical motions, which
belong to IMSMs under scrutiny. When representing them in the form
of the Sylvester conditions, we have:
\begin{gather}
 1) \ \  \frac{\alpha^2 \lambda_0^2 + \lambda_1^2}{2 \lambda_0} > 0,
\nonumber \\
 2) \ \ \frac{(\alpha^2 \lambda_0^2 + \lambda_1^2)
  \big(\lambda_1^4 + \alpha^2 \lambda_0 (\alpha^2 \lambda_0^2 + \lambda_1^2)
  \big(\lambda_0 - \lambda_3 {s_3^0}^2\big)\big)} {4 \lambda_0^2 \lambda_1^2} > 0, \nonumber \\
 3) \ \
  \frac{1} {8 \lambda_0^2 \lambda_1^2}
  \big(\alpha^2 \lambda_0^3 - (\alpha^2 \lambda_0^2 + \lambda_1^2) \lambda_3 {s_3^0}^2\big)
  \big(\lambda_1^4 + \alpha^2 \lambda_0 (\alpha^2 \lambda_0^2 + \lambda_1^2)
  \big(\lambda_0 - \lambda_3 {s_3^0}^2\big)\big) > 0, \nonumber \\
\nonumber \\
 4) \ \  \frac{(\alpha^2 \lambda_0^2 + \lambda_1^2)^3
    \big(\lambda_1^2 \lambda_3 {s_3^0}^2
    + \alpha^2 \lambda_0^2 \big(\lambda_3 {s_3^0}^2 - \lambda_0\big)\big)^2}
    {16 \alpha^2 \lambda_0^4 \lambda_1^4} > 0.\label{15}
\end{gather}
A standard software package ``Algebra InequalitySolve'' of CAS
``Mathematica'' was used for the purpose of verif\/ication of
compatibility for this system of inequalities. Its application
to~(\ref{15}) showed that the inequalities are compatible when:
\begin{gather*}
 \alpha < 0 \vee \lambda_0 > 0 \\
 \vee \Bigg ( \lambda_1 < 0 \vee \Bigg(\lambda_3 < 0 \wedge \lambda_3 > 0 \vee
    -\frac{\alpha \lambda_0^{3/2}}
    {\sqrt{(\alpha^2 \lambda_0^2 + \lambda_1^2) \lambda_3}} <
    s_3^0 < \frac{\alpha \lambda_0^{3/2}}{\sqrt{(\alpha^2 \lambda_0^2
    + \lambda_1^2) \lambda_3}} \Bigg)
  \\
  \wedge \lambda_1 > 0
    \vee \Bigg(\lambda_3 < 0 \wedge \lambda_3 > 0 \vee -\frac{\alpha \lambda_0^{3/2}}
    {\sqrt{(\alpha^2 \lambda_0^2 + \lambda_1^2) \lambda_3}} <
    s_3^0 < \frac{\alpha \lambda_0^{3/2}}{\sqrt{(\alpha^2 \lambda_0^2 +
    \lambda_1^2) \lambda_3}} \Bigg) \Bigg )
    \\
 \wedge \alpha > 0 \vee  \lambda_0 > 0 \\
 \vee \Bigg (\lambda_1 < 0 \vee \Bigg(\lambda_3 < 0 \vee \lambda_3 > 0 \wedge
     - \frac{\alpha \lambda_0^{3/2}}
        {\sqrt{(\alpha^2 \lambda_0^2 + \lambda_1^2) \lambda_3}} <
       s_3^0 < \frac{\alpha \lambda_0^{3/2}}{\sqrt{(\alpha^2 \lambda_0^2 +
           \lambda_1^2) \lambda_3}} \Bigg) \\
  \wedge  \lambda_1 > 0 \vee
     \Bigg( \lambda_3 < 0 \vee \lambda_3 > 0 \wedge -\frac{\alpha \lambda_0^{3/2}}
        {\sqrt{(\alpha^2 \lambda_0^2 + \lambda_1^2) \lambda_3}} <
       s_3^0 < \frac{\alpha \lambda_0^{3/2}}{\sqrt{(\alpha^2 \lambda_0^2 +
           \lambda_1^2) \lambda_3}} \Bigg) \Bigg ).
\end{gather*}
 Hence, stable in the sense of Lyapunov are only those
body's helical motions for which the parameter $s_3^0$ satisf\/ies
the latter conditions.

Comparison of stability conditions for the zero solution with the
stability conditions for IMSMs (\ref{117}) gives evidence that
satisfaction of the former conditions implies satisfaction of the
latter. Consequently, the stable IMSMs adjunct to the stable zero
solution.

We also investigated stability of IMSMs (\ref{1172}). We proved
that some of these IMSMs are unstable with respect to the f\/irst
approximation. Investigation of stability for other solutions
necessitates involving of higher-order terms in the integral $K$
 expansion. But this problem was not considered here.

\section{Conclusion}
 The paper  represented some results of qualitative analysis of the
dif\/ferential equations describing the motion of a rigid body in
ideal f\/luid. The Kirchhof\/f's dif\/ferential equations were
considered in Sokolov's case. In this case, these have 4 algebraic
f\/irst integrals (three quadratic ones and one 4th degree
integral) and represent a completely integrable system.

 Routh--Lyapunov's method  was used to analyze the set of solutions
 of the equations. The method proposes a technique of f\/inding both
 the stationary solutions and the invariant manifolds of motion
 equations, when the equations have a suf\/f\/iciently large number of f\/irst
 integrals. The stationary solutions and invariant manifolds, which were
 obtained by this method, may be investigated for stability by the 2nd Lyapunov's
 method. Furthermore, the corresponding f\/irst integrals can be used here
 as Lyapunov functions.

 The paper considers a rather typical case when the problem of f\/inding
stationary solutions is reduced to solving a nonlinear system of
equations. Some interesting cases of stability investigation for
conservative systems are also given. For example, the property of
asymptotic stability of the equilibrium state for a vector f\/ield
on the 1-dimensional invariant manifold of a conservative system
was used for proving instability of the stationary solution (which
corresponds to this equilibrium state) within the whole system's
phase space.
 Some examples of branching invariant manifolds of various dimensions
 and investigation of their stability were considered.

 The results of this work were obtained with the use of computer algebra
tools. Such investigations cannot probably be conducted within an
acceptable time without them. To ground this statement, we gave a
Gr\"{o}bner basis in Appendix C. This basis was constructed and
used for f\/inding solutions of a nonlinear system of algebraic
equations arising in
 computations.

 The results of qualitative analysis of Kirchhof\/f's dif\/ferential equations
represented in the paper give evidence that the technique of
investigation of mechanical systems, which is based on a
combination of classical methods of rigid body dynamics and
computer algebra methods, is rather ef\/f\/icient and may be used
for investigations of above type problems.

\appendix

\section{Appendix}
The following concepts were used in the paper.

\begin{definition}
The solutions of dif\/ferential equations, on which the f\/irst integral (an
element of the algebra of problem's f\/irst integrals) assumes a stationary value, are called
the stationary solutions.
\end{definition}

This means the following. The stationary solution satisf\/ies the equations
which are the result of equating all the partial derivatives of
the f\/irst integral with respect to the problem's variables to
zero (stationary conditions). Such solutions simultaneously are
the solutions of initial dif\/ferential equations. This follows
from the Lyapunov's theorem \cite{Lyapunov&1956}.

\begin{definition}
Manifolds, whose equations satisfy the stationary conditions for
some f\/irst integral, are called the invariant manifolds of
steady motions (IMSMs).
\end{definition}
The proof of invariance of such manifolds with respect to initial
dif\/ferential equations can be found in Appendix B.
 Note, that manifolds are understood as the sets of dimension larger
than zero.

 If some dif\/ferential equations have an IMSM then by these equations it is
possible to def\/ine a~vector f\/ield on this IMSM. This procedure
is called a reduction of the initial system.  The vector f\/ield,
in turn, can have f\/irst integrals which may be used for
f\/inding stationary solutions and the 2nd-level IMSM (on the
given IMSM).

 Stationary solutions and manifolds are suitable in stability investigations
by Lyapunov's second method, because the expansions of the
corresponding f\/irst integrals in the neighbourhood of these
solutions and manifolds in Taylor series do not contain linear
terms.

\section{Appendix}
 More exactly, we use the following theorem:
\begin{theorem}
\label{theorem1}
 If partial derivatives of the first integral $ V(x,t) $ of the system $ \dot{x}_{i}=
X_{i}(x,t) \ (i=1,\ldots,n)$ with respect to the problem's phase
variables have the form
\[
\frac{ \partial V}{ \partial x_{i}}= \sum_{l=1}^{k}a_{il}(x,t)
\varphi_{l}(x,t) + \sum_{l=1}^{k} \sum_{p=1}^{k}a_{il}(x,t)
\varphi_{l}(x,t) \varphi_{p}(x,t)+\cdots, \qquad i=1,\ldots,n,
\]
and the rank of the matrix $ \Vert a_{il}(x,t) \Vert$ is ``$k$''
on the manifold
 $ \varphi_{l}(x,t)=0$,  $l=1,\ldots,k$, then the manifold $ \varphi_{l}(x,t)=0$,
$l=1,\ldots,k$, is invariant for the initial system of
differential
 equations.
\end{theorem}
\begin{proof}
Let the system of dif\/ferential equations
\[
 \dot{x_{i}}=X_{i}(x,t), \qquad i=1,2,\ldots, n,
\]
  have the f\/irst integral $ V(x,t) $ and
\[
\frac{ \partial V}{ \partial x_{i}}= \sum_{l=1}^{k}a_{il}(x,t)
\varphi_{l}(x,t) + \sum_{l=1}^{k} \sum_{p=1}^{k}a_{ilp}(x,t)
\varphi_{l}(x,t) \varphi_{p}(x,t)+\cdots,\qquad i=1,\ldots,n.
\]
Since $ V(x,t) $ is the f\/irst integral, we have
\[
\frac {dV}{dt}= \sum_{j=1}^{n} \frac { \partial V}{ \partial
x_{j}}X_{j} + \frac { \partial V}{ \partial t}=0.
\]
When dif\/ferentiating the latter identity with respect to $ x_{i}
$, we obtain the system of equations
\[
\sum_{j=1}^{n} \frac { \partial}{ \partial x_{i}}\left( \frac {
\partial V}
 { \partial x_{j}}\right)X_{j} + \frac { \partial}{ \partial x_{i}}
\frac { \partial V}{ \partial t} + \sum_{j=1}^{n} \frac { \partial
V}
 { \partial x_{j}}\frac { \partial X_{j}}{ \partial x_{i}}=0,\qquad i=1,\ldots,n.
\]
Having changed the order of dif\/ferentiation in the above system,
we have
\[
\sum_{j=1}^{n} \frac { \partial}{ \partial x_{j}}\left( \frac {
\partial V}
 { \partial x_{i}}\right)X_{j} + \frac { \partial}{ \partial t}
\frac { \partial V}{ \partial x_{i}} =- \sum_{j=1}^{n} \frac {
\partial V}
 { \partial x_{j}} \frac { \partial X_{j}}{ \partial x_{i}},
 \qquad i=1,\ldots,n.
\]
Now substitute the expression for ${\partial V}/{ \partial x_{i}}$
into the latter formulae and perform dif\/ferentiation with
respect to $x_i$. After trivial transformations we have:
\[
\sum_{l=1}^{k} a_{il}\left( \sum_{j=1}^{n} \frac { \partial
\varphi_{l}}
 { \partial x_{j}}X_{j} + \frac { \partial \varphi_{l}}{ \partial t}\right)=
 F_{i}( \varphi_{l}), \qquad i=1,\ldots,n,
\]
where $ F_{i}(0)=0$.

  When $ \varphi_{p} =0$, $p=1,\ldots,k$, the latter system transforms
  into the system of $n$ linear homogeneous equations:
\[
 \sum_{l=1}^{k} a_{il}\left( \sum_{j=1}^{n} \frac { \partial \varphi_{l}}
 { \partial x_{j}}X_{j} + \frac { \partial \varphi_{l}}{ \partial t}\right)= 0
\qquad i=1,\ldots,n.
\]
 Under the condition that the rank of the matrix $\Vert a_{il}(x,t) \Vert$
is $ k $ on the manifold $ \varphi_{p} =0$, $p=1,\ldots,k$, this
system has only the following trivial solution:
\[
 \sum_{j=1}^{n} \frac { \partial \varphi_{l}}
 { \partial x_{j}} X_{j} + \frac { \partial \varphi_{l}}{ \partial t}= 0,
\qquad l=1,\ldots, k.
\]
The latter proves that the manifold $ \varphi_{p} = 0$,
$p=1,\ldots,k$ is invariant for the initial system of
dif\/ferential equations.
\end{proof}

\section{Appendix}
The Gr\"{o}bner basis constructed for the system (\ref{0404}):
\begin{gather*}
  s_2 s_3 \big \{2 \lambda_2 (\alpha^4 \lambda_0^3 + 4
\alpha^2 \lambda_0^2 \lambda_2 +
       2 \lambda_2 (\lambda_1^2 + 2 \lambda_0 \lambda_2)) +
     4 \alpha \lambda_1 \lambda_2 (\alpha^2 \lambda_0
     + 2 \lambda_2) \lambda_3 s_1 s_3 \\
\qquad {}- 2 \lambda_2 (\alpha^4 \lambda_0^2
     + 4 \lambda_2^2
  - 2 \alpha^2 (\lambda_1^2 - 2 \lambda_0 \lambda_2)) \lambda_3
s_3^2 +
     \alpha^2 (\alpha^2 \lambda_0 + 2 \lambda_2)^2 \lambda_3^2 s_1^2 s_3^2\\
\qquad {} +
     2 \alpha^3 \lambda_1 (\alpha^2 \lambda_0 + 2 \lambda_2) \lambda_3^2 s_1 s_3^3 +
     \alpha^4 \lambda_1^2 \lambda_3^2 s_3^4 \big \} = 0, \\
  ((\alpha^2 \lambda_0 + 2 \lambda_2) s_1 + \alpha \lambda_1
 s_3) \big \{2 \lambda_2 (\alpha^4 \lambda_0^3 + 4 \alpha^2 \lambda_0^2 \lambda_2 +
      2 \lambda_2 (\lambda_1^2 + 2 \lambda_0 \lambda_2))\\
 \qquad{} +
    4 \alpha \lambda_1 \lambda_2 (\alpha^2 \lambda_0 + 2 \lambda_2) \lambda_3 s_1 s_3
 - 2 \lambda_2 (\alpha^4 \lambda_0^2 + 4 \lambda_2^2 -
      2 \alpha^2 (\lambda_1^2 - 2 \lambda_0 \lambda_2)) \lambda_3 s_3^2 \\
 \qquad{}+
    \alpha^2 (\alpha^2 \lambda_0 + 2 \lambda_2)^2 \lambda_3^2 s_1^2 s_3^2 +
    2 \alpha^3 \lambda_1 (\alpha^2 \lambda_0 + 2 \lambda_2) \lambda_3^2 s_1 s_3^3
    +  \alpha^4 \lambda_1^2 \lambda_3^2 s_3^4 \big \} = 0, \\
  s_2 \big \{2 \lambda_2 (\alpha^2 \lambda_0 + 2 \lambda_2)^3
      (\lambda_1^2 + 2 \lambda_0 \lambda_2) + 4 \alpha^2 \lambda_0 \lambda_2
      (\alpha^2 \lambda_0 + 2 \lambda_2) (\alpha^2 \lambda_0 - \alpha \lambda_1 +
       2 \lambda_2)\\
 \qquad{}\times  (\alpha (\alpha \lambda_0 + \lambda_1)
 + 2 \lambda_2) \lambda_3 s_1^2- \alpha^4 \lambda_0 (\alpha^2
\lambda_0 + 2 \lambda_2)^3 \lambda_3^2 s_1^4 +
     4 \alpha^2 \lambda_0 \lambda_2 (\alpha^2 \lambda_0 + 2 \lambda_2)^3 \lambda_3 s_2^2 \\
\qquad{}+
     \alpha^4 \lambda_0 (\alpha^2 \lambda_0 + 2 \lambda_2)^3 \lambda_3^2 s_2^4
 + 4 \alpha \lambda_1 \lambda_2 (2 \alpha^6 \lambda_0^3 + 16
\alpha^2 \lambda_0
        \lambda_2^2 + 8 \lambda_2^3  \\
\qquad{}+ \alpha^4 \lambda_0 (10 \lambda_0 \lambda_2
        -\lambda_1^2)) \lambda_3 s_1 s_3- 4 \alpha^5 \lambda_0 \lambda_1
      (\alpha^2 \lambda_0 + 2 \lambda_2)^2 \lambda_3^2 s_1^3 s_3
 - 2 \lambda_2 (\alpha^2 \lambda_0 + 2 \lambda_2)\\
 \qquad{}\times (8
\lambda_2^3 -
       4 \alpha^2 \lambda_2 (\lambda_1^2 - 2 \lambda_0 \lambda_2) +
       \alpha^4 \lambda_0 (2 \lambda_0 \lambda_2-\lambda_1^2)) \lambda_3 s_3^2 +
     \alpha^2 (\alpha^2 \lambda_0 + 2 \lambda_2)\\
\qquad{}\times (8 \alpha^2 \lambda_0
 \times \lambda_2^2 + 8 \lambda_2^3 + \alpha^4 \lambda_0 (-5
\lambda_1^2 + 2 \lambda_0 \lambda_2))
      \lambda_3^2 s_1^2 s_3^2 + 2 \alpha^3 \lambda_1 (8 \alpha^2 \lambda_0 \lambda_2^2 +
       8 \lambda_2^3 \\
\qquad{}\times+ \alpha^4 \lambda_0 (2 \lambda_0
\lambda_2-\lambda_1^2))
      \lambda_3^2 s_1 s_3^3  + 2 \alpha^4 \lambda_1^2 \lambda_2
      (\alpha^2 \lambda_0 + 2 \lambda_2) \lambda_3^2 s_3^4 \big \} = 0, \\
 s_3 \big \{2 \lambda_2 (\alpha^4 \lambda_0^3 + 4 \alpha^2
\lambda_0^2 \lambda_2 +
     2 \lambda_2 (\lambda_1^2 + 2 \lambda_0 \lambda_2)) +
   4 \alpha \lambda_1 \lambda_2 (\alpha^2 \lambda_0 + 2 \lambda_2) \lambda_3 s_1 s_3\\
\qquad{}   - 2 \lambda_2 (\alpha^4 \lambda_0^2 + 4 \lambda_2^2
 -2 \alpha^2 (\lambda_1^2 - 2 \lambda_0 \lambda_2)) \lambda_3
s_3^2 +
   \alpha^2 (\alpha^2 \lambda_0 + 2 \lambda_2)^2 \lambda_3^2 s_1^2 s_3^2\\
\qquad{} +
   2 \alpha^3 \lambda_1 (\alpha^2 \lambda_0 + 2 \lambda_2) \lambda_3^2 s_1 s_3^3 +
   \alpha^4 \lambda_1^2 \lambda_3^2 s_3^4 \big \}
\big \{(\alpha^2 \lambda_0 + 2 \lambda_2)^4 (\alpha^2
\lambda_0^2 + \lambda_1^2 +     4 \lambda_0 \lambda_2) \\
\qquad{}- 4 \lambda_2 (\alpha^2 \lambda_0 + 2 \lambda_2)^2
    (\alpha^4 \lambda_0^2 + 4 \lambda_2^2 +
     \alpha^2 (\lambda_1^2 + 4 \lambda_0 \lambda_2)) \lambda_3
     s_3^2  - \alpha^2 (\alpha^4 \lambda_0^2 + 4 \lambda_2^2\\
\qquad{}+
      \alpha^2 (\lambda_1^2 + 4 \lambda_0 \lambda_2))^2 \lambda_3^2
      s_3^4 \big \}     = 0, \\
  2 \lambda_1 \lambda_2 (\alpha^2 \lambda_0 + 2 \lambda_2)^2
r_2 - 2 \lambda_0 \lambda_2 (\alpha^2 \lambda_0 + 2 \lambda_2)^2
s_2 -
   \lambda_0 (\alpha^3 \lambda_0 + 2 \alpha \lambda_2)^2 \lambda_3 s_1^2 s_2\\
\qquad{}-
   \lambda_0 (\alpha^3 \lambda_0 + 2 \alpha \lambda_2)^2 \lambda_3   s_2^3
   - \alpha \lambda_1 (\alpha^2 \lambda_0 + 2
\lambda_2)^2 \lambda_3 s_1 s_2 s_3 +
   2 \lambda_2 (\alpha^2 \lambda_0 - \alpha \lambda_1 + 2 \lambda_2)\\
\qquad{}\times
    (\alpha (\alpha \lambda_0 + \lambda_1) + 2 \lambda_2) \lambda_3 s_2 s_3^2 -
   \alpha^2 (\alpha^2 \lambda_0  + 2 \lambda_2)^2 \lambda_3^2 s_1^2 s_2 s_3^2\\
\qquad{} -
   2 \alpha^3 \lambda_1 (\alpha^2 \lambda_0 + 2 \lambda_2) \lambda_3^2 s_1 s_2 s_3^3 -
   \alpha^4 \lambda_1^2 \lambda_3^2 s_2 s_3^4 = 0, \\
  s_3 \big \{-2 \lambda_2 (\alpha^2 \lambda_0 + 2
\lambda_2)^8 (\alpha^2 \lambda_0^2 +
      \lambda_1^2 + 4 \lambda_0 \lambda_2) + (\alpha^2 \lambda_0 + 2 \lambda_2)^2
     (\alpha^{16} \lambda_0^8 + 256 \lambda_2^8 \\
\qquad{}+ 16 \alpha^{12} \lambda_0^5 \lambda_2
       (\lambda_1^2 + 7 \lambda_0  \lambda_2) + \alpha^{14} \lambda_0^6
       (\lambda_1^2 + 16 \lambda_0 \lambda_2) + 64 \alpha^2 \lambda_2^6
       (3 \lambda_1^2 + 16 \lambda_0 \lambda_2)\\
\qquad{} + 32 \alpha^4 \lambda_2^4
       (\lambda_1^4 + 16 \lambda_0 \lambda_1^2 \lambda_2 + 56 \lambda_0^2 \lambda_2^2) +
      4 \alpha^{10} \lambda_0^3 \lambda_2  (\lambda_1^4 + 25 \lambda_0 \lambda_1^2
         \lambda_2 + 112 \lambda_0^2 \lambda_2^2)\\
\qquad{} + 8 \alpha^8 \lambda_0^2 \lambda_2^2
       (3 \lambda_1^4 + 40 \lambda_0 \lambda_1^2 \lambda_2
  + 140 \lambda_0^2 \lambda_2^2) + 4 \alpha^6 \lambda_2^2
       (\lambda_1^6 + 12 \lambda_0 \lambda_1^4 \lambda_2
 + 140 \lambda_0^2   \lambda_1^2 \lambda_2^2\\
 \qquad{} + 448 \lambda_0^3 \lambda_2^3))
\lambda_3 s_1^2 +
    (\alpha^2 \lambda_0 + 2 \lambda_2)^2 (\alpha^{16} \lambda_0^8 + 256 \lambda_2^8 +
      16 \alpha^{12} \lambda_0^5 \lambda_2 (\lambda_1^2 + 7 \lambda_0 \lambda_2)\\
\qquad{} +
      \alpha^{14} \lambda_0^6 (\lambda_1^2 + 16 \lambda_0
      \lambda_2)+ 64 \alpha^2 \lambda_2^6 (3 \lambda_1^2 + 16
\lambda_0 \lambda_2) +
      32 \alpha^4 \lambda_2^4 (\lambda_1^4 + 16 \lambda_0 \lambda_1^2 \lambda_2 +
        56 \lambda_0^2 \lambda_2^2)\\
\qquad{} + 4 \alpha^{10} \lambda_0^3 \lambda_2
       (\lambda_1^4 + 25 \lambda_0 \lambda_1^2 \lambda_2
     + 112 \lambda_0^2  \lambda_2^2)+ 8 \alpha^8 \lambda_0^2 \lambda_2^2 (3
\lambda_1^4 + 40 \lambda_0 \lambda_1^2
         \lambda_2 + 140 \lambda_0^2 \lambda_2^2)\\
\qquad{} + 4 \alpha^6 \lambda_2^2
       (\lambda_1^6 + 12 \lambda_0 \lambda_1^4 \lambda_2 + 140 \lambda_0^2 \lambda_1^2
         \lambda_2^2 + 448 \lambda_0^3 \lambda_2^3)) \lambda_3 s_2^2
 + 2 \alpha \lambda_1 (\alpha^2 \lambda_0 + 2 \lambda_2)^3\\
\qquad{}\times
     (\alpha^{12} \lambda_0^6 - 96 \alpha^2 \lambda_0 \lambda_2^5 - 64 \lambda_2^6 -
      16 \alpha^4 \lambda_0 \lambda_2^3 (\lambda_0 \lambda_2-\lambda_1^2) +
      2 \alpha^6 \lambda_0 \lambda_2 (\lambda_1^2 + 4 \lambda_0
      \lambda_2)\\
\qquad{}\times (\lambda_1^2  + 6 \lambda_0 \lambda_2) + 4 \alpha^8
\lambda_0^3 \lambda_2
       (2 \lambda_1^2 + 9 \lambda_0 \lambda_2) + \alpha^{10} \lambda_0^4
       (\lambda_1^2 + 10 \lambda_0 \lambda_2)) \lambda_3 s_1 s_3\\
\qquad{} +
    (\alpha^{18} \lambda_0^8 \lambda_1^2 + 2048 \lambda_2^10
    - 256 \alpha^2 \lambda_2^8 (\lambda_1^2  - 32 \lambda_0 \lambda_2) +
      \alpha^{16} \lambda_0^6 (\lambda_1^4 + 12 \lambda_0 \lambda_1^2 \lambda_2 +
        8 \lambda_0^2 \lambda_2^2)\\
\qquad{} + 8 \alpha^{14} \lambda_0^5 \lambda_2
       (\lambda_1^4 + 7 \lambda_0 \lambda_1^2 \lambda_2
       + 16 \lambda_0^2 \lambda_2^2) + 64 \alpha^4 \lambda_2^6 (
        224 \lambda_0^2 \lambda_2^2 -3 \lambda_1^4 - 12 \lambda_0 \lambda_1^2 \lambda_2 )\\
\qquad{} + 4 \alpha^{12} \lambda_0^4 \lambda_2^2
       (3 \lambda_1^4 + 28 \lambda_0 \lambda_1^2 \lambda_2 +
        224 \lambda_0^2 \lambda_2^2) - 32 \alpha^6 \lambda_2^4
       (\lambda_1^6 + 12 \lambda_0 \lambda_1^4 \lambda_2 + 28 \lambda_0^2
\lambda_1^2
         \lambda_2^2 \\
\qquad{}  - 448 \lambda_0^3 \lambda_2^3) + 8 \alpha^{10}
\lambda_0^2
       \lambda_2^2 (
        448 \lambda_0^3 \lambda_2^3 -\lambda_1^6 - 8 \lambda_0 \lambda_1^4 \lambda_2) - 4 \alpha^8 \lambda_2^2
       (\lambda_1^8 + 8 \lambda_0 \lambda_1^6 \lambda_2 + 68 \lambda_0^2 \lambda_1^4
         \lambda_2^2\\
\qquad{}+ 112 \lambda_0^3 \lambda_1^2 \lambda_2^3
  - 2240 \lambda_0^4 \lambda_2^4)) \lambda_3 s_3^2 -
    \alpha^2 (\alpha^2 \lambda_0 + 2 \lambda_2)^2 (\alpha^4 \lambda_0^2 +
      4 \lambda_2^2 + \alpha^2 (\lambda_1^2 + 4 \lambda_0 \lambda_2))\\
\qquad{}\times
     (\alpha^{10} \lambda_0^5 + 64 \lambda_2^5 + \alpha^8 \lambda_0^3
       (\lambda_1^2 + 12 \lambda_0 \lambda_2) + 8 \alpha^2 \lambda_2^3
       (\lambda_1^2 + 18 \lambda_0 \lambda_2) + 2 \alpha^6 \lambda_0^2 \lambda_2
       (3 \lambda_1^2  \\
\qquad{}+ 28 \lambda_0 \lambda_2)+ 2 \alpha^4 \lambda_2
       (\lambda_1^4 + 6 \lambda_0 \lambda_1^2 \lambda_2
       + 64 \lambda_0^2 \lambda_2^2))
     \lambda_3^2 s_1^2 s_3^2  - 2 \alpha^3 \lambda_1 (\alpha^2 \lambda_0 + 2 \lambda_2)
     (\alpha^4 \lambda_0^2 + 4 \lambda_2^2\! \\
\qquad{}+
      \alpha^2 (\lambda_1^2 + 4 \lambda_0 \lambda_2)) (\alpha^{10} \lambda_0^5 +
      96 \lambda_2^5 + 8 \alpha^6 \lambda_0^2 \lambda_2 (\lambda_1^2
      + 9 \lambda_0 \lambda_2) + 16 \alpha^2 \lambda_2^3
  (\lambda_1^2+ 13 \lambda_0 \lambda_2)\\
\qquad{} + \alpha^8 \lambda_0^3
       (\lambda_1^2 + 14 \lambda_0 \lambda_2) + 2 \alpha^4 \lambda_2
       (\lambda_1^4 + 10 \lambda_0 \lambda_1^2 \lambda_2 + 88 \lambda_0^2 \lambda_2^2))
     \lambda_3^2 s_1 s_3^3 + \alpha^2 (\alpha^4 \lambda_0^2 + 4 \lambda_2^2\\
\qquad{} +      \alpha^2  (\lambda_1^2+ 4 \lambda_0 \lambda_2))
(128 \lambda_2^7 - 96 \alpha^2 \lambda_2^5 (\lambda_1^2 - 4
\lambda_0 \lambda_2) +
      \alpha^{12} \lambda_0^5 (2 \lambda_0 \lambda_2 - \lambda_1^2)\\
\qquad{}
 + \alpha^{10} \lambda_0^3 (
        24 \lambda_0^2 \lambda_2^2 - \lambda_1^4
        - 14 \lambda_0 \lambda_1^2 \lambda_2 )     + 2 \alpha^8 \lambda_0^2 \lambda_2
       (
        60 \lambda_0^2 \lambda_2^2 - 5 \lambda_1^4
        - 36 \lambda_0 \lambda_1^2 \lambda_2 )\\
 \qquad{} + 8 \alpha^4 \lambda_2^3
       ( 60 \lambda_0^2 \lambda_2^2
        -3 \lambda_1^4 - 26 \lambda_0 \lambda_1^2 \lambda_2) - 2
\alpha^6 \lambda_2
       (\lambda_1^6 + 14 \lambda_0 \lambda_1^4 \lambda_2 + 88 \lambda_0^2
 \lambda_1^2
         \lambda_2^2\\
 \qquad{} - 160 \lambda_0^3 \lambda_2^3)) \lambda_3^2 s_3^4 -
    (\alpha^3 \lambda_0 + 2 \alpha \lambda_2)^4
    (\alpha^4 \lambda_0^2 + 4 \lambda_2^2 + \alpha^2 (\lambda_1^2
+  4 \lambda_0 \lambda_2))^2
 \lambda_3^3 s_1^2 s_3^4\\
 \qquad{} -
    2 \alpha^5 \lambda_1 (\alpha^2 \lambda_0 + 2 \lambda_2)^3
     (\alpha^4 \lambda_0^2 + 4 \lambda_2^2 + \alpha^2 (\lambda_1^2 +
         4 \lambda_0 \lambda_2))^2 \lambda_3^3
         s_1 s_3^5 - \alpha^6 \lambda_1^2 (\alpha^2 \lambda_0 + 2
\lambda_2)^2\\
\qquad{}\times
     (\alpha^4 \lambda_0^2 + 4 \lambda_2^2 + \alpha^2 (\lambda_1^2
  +4 \lambda_0 \lambda_2))^2 \lambda_3^3 s_3^6 \big \} = 0,
\\
  -2 \lambda_1^2 \lambda_2 (\alpha^2 \lambda_0 + 2 \lambda_2)
    (\alpha^{16} \lambda_0^8 + 256 \lambda_2^8 + 16 \alpha^{12} \lambda_0^5 \lambda_2
      (\lambda_1^2 + 7 \lambda_0 \lambda_2) + \alpha^{14} \lambda_0^6
      (\lambda_1^2 + 16 \lambda_0 \lambda_2)\\
\qquad{} + 64 \alpha^2 \lambda_2^6
      (3 \lambda_1^2  + 16 \lambda_0 \lambda_2) + 32 \alpha^4 \lambda_2^4
      (\lambda_1^4 + 16 \lambda_0 \lambda_1^2 \lambda_2 + 56 \lambda_0^2 \lambda_2^2)\\
\qquad{} +
     4 \alpha^{10} \lambda_0^3 \lambda_2 (\lambda_1^4 + 25 \lambda_0 \lambda_1^2 \lambda_2 +
       112 \lambda_0^2 \lambda_2^2) + 8 \alpha^8 \lambda_0^2 \lambda_2^2
      (3 \lambda_1^4  + 40 \lambda_0 \lambda_1^2 \lambda_2 +
       140 \lambda_0^2 \lambda_2^2)\\
\qquad{} + 4 \alpha^6 \lambda_2^2
      (\lambda_1^6 + 12 \lambda_0 \lambda_1^4 \lambda_2 + 140 \lambda_0^2 \lambda_1^2
        \lambda_2^2 + 448 \lambda_0^3 \lambda_2^3)) s_1 -
   (\alpha^2 \lambda_0 + 2 \lambda_2)^3 (\alpha^{16} \lambda_0^8
 + 256 \lambda_2^8\\
 \qquad{} +
     16 \alpha^{12} \lambda_0^5 \lambda_2 (\lambda_1^2 + 7 \lambda_0 \lambda_2) +
     \alpha^{14} \lambda_0^6 (\lambda_1^2 + 16 \lambda_0 \lambda_2) +
     64 \alpha^2 \lambda_2^6 (3 \lambda_1^2 + 16 \lambda_0 \lambda_2)\\
 \qquad{} +     32 \alpha^4 \lambda_2^4 (\lambda_1^4
 + 16 \lambda_0 \lambda_1^2 \lambda_2 +
       56 \lambda_0^2 \lambda_2^2) + 4 \alpha^{10} \lambda_0^3 \lambda_2
      (\lambda_1^4 + 25 \lambda_0 \lambda_1^2 \lambda_2 + 112 \lambda_0^2 \lambda_2^2)\\
 \qquad{} +
     8 \alpha^8 \lambda_0^2 \lambda_2^2 (3 \lambda_1^4 + 40 \lambda_0 \lambda_1^2
        \lambda_2 + 140 \lambda_0^2 \lambda_2^2)  + 4 \alpha^6 \lambda_2^2
      (\lambda_1^6 + 12 \lambda_0 \lambda_1^4 \lambda_2 + 140 \lambda_0^2 \lambda_1^2
        \lambda_2^2\\
 \qquad{} + 448 \lambda_0^3 \lambda_2^3)) \lambda_3 s_1^3 -
   (\alpha^2 \lambda_0 + 2 \lambda_2)^3 (\alpha^{16} \lambda_0^8 + 256 \lambda_2^8 +
     16 \alpha^{12} \lambda_0^5 \lambda_2
 (\lambda_1^2 + 7 \lambda_0 \lambda_2)\\
 \qquad{} +
     \alpha^{14} \lambda_0^6 (\lambda_1^2 + 16 \lambda_0 \lambda_2) +
     64 \alpha^2 \lambda_2^6 (3 \lambda_1^2 + 16 \lambda_0 \lambda_2) +
     32 \alpha^4 \lambda_2^4 (\lambda_1^4 + 16 \lambda_0 \lambda_1^2 \lambda_2 +
       56 \lambda_0^2 \lambda_2^2) \\
\qquad{} + 4 \alpha^{10} \lambda_0^3 \lambda_2
      (\lambda_1^4 + 25 \lambda_0 \lambda_1^2 \lambda_2 + 112 \lambda_0^2 \lambda_2^2) +
     8 \alpha^8 \lambda_0^2 \lambda_2^2 (3 \lambda_1^4 + 40 \lambda_0 \lambda_1^2
        \lambda_2 + 140 \lambda_0^2 \lambda_2^2)\\
\qquad{} + 4 \alpha^6 \lambda_2^2
      (\lambda_1^6 + 12 \lambda_0 \lambda_1^4
  \times \lambda_2 + 140 \lambda_0^2 \lambda_1^2
        \lambda_2^2 + 448 \lambda_0^3 \lambda_2^3)) \lambda_3 s_1 s_2^2 -
   2 \alpha \lambda_1 \lambda_2 (\alpha^4 \lambda_0^3 \\
  \qquad{}+ 4 \alpha^2 \lambda_0^2 \lambda_2 +
     2 \lambda_2 (\lambda_1^2 + 2 \lambda_0 \lambda_2)) (\alpha^4 \lambda_0^2
     + 4 \lambda_2^2 + \alpha^2  (\lambda_1^2 + 4 \lambda_0 \lambda_2))
    (\alpha^{10} \lambda_0^5 + 64 \lambda_2^5\\
  \qquad{} + \alpha^8 \lambda_0^3
      (\lambda_1^2 + 12 \lambda_0 \lambda_2) + 8 \alpha^2 \lambda_2^3
      (\lambda_1^2 + 18 \lambda_0 \lambda_2) + 2 \alpha^6 \lambda_0^2 \lambda_2
      (3 \lambda_1^2 + 28 \lambda_0 \lambda_2) \\
 \qquad{}+ 2 \alpha^4 \lambda_2
      (\lambda_1^4  + 6 \lambda_0 \lambda_1^2 \lambda_2 + 64 \lambda_0^2
\lambda_2^2))
    s_3 - 2 \alpha \lambda_1 (\alpha^2 \lambda_0 + 2 \lambda_2)^2
    (\alpha^{16} \lambda_0^8 + 256 \lambda_2^8 \\
\qquad{}+ 16 \alpha^{12} \lambda_0^5 \lambda_2
      (\lambda_1^2 + 7 \lambda_0 \lambda_2)  + \alpha^{14} \lambda_0^6
      (\lambda_1^2 + 16 \lambda_0 \lambda_2) + 64 \alpha^2 \lambda_2^6
      (3 \lambda_1^2 + 16 \lambda_0 \lambda_2)\\
\qquad{} + 32 \alpha^4 \lambda_2^4
      (\lambda_1^4 + 16 \lambda_0 \lambda_1^2 \lambda_2
      + 56 \lambda_0^2 \lambda_2^2) + 4 \alpha^{10} \lambda_0^3 \lambda_2 (\lambda_1^4
   + 25 \lambda_0 \lambda_1^2 \lambda_2 +
       112 \lambda_0^2 \lambda_2^2)\\
\qquad{} + 8 \alpha^8 \lambda_0^2 \lambda_2^2
      (3 \lambda_1^4 + 40 \lambda_0 \lambda_1^2 \lambda_2 +
       140 \lambda_0^2 \lambda_2^2) + 4 \alpha^6 \lambda_2^2
      (\lambda_1^6 + 12 \lambda_0 \lambda_1^4 \lambda_2 + 140 \lambda_0^2 \lambda_1^2
        \lambda_2^2 \\
\qquad{}  + 448 \lambda_0^3 \lambda_2^3)) \lambda_3 s_1^2 s_3 -
   \alpha^2 \lambda_1^2 (\alpha^2 \lambda_0 + 2 \lambda_2)
    (\alpha^{16} \lambda_0^8 + 1280 \lambda_2^8 + \alpha^{14} \lambda_0^6
      (\lambda_1^2 + 20 \lambda_0 \lambda_2)\\
\qquad{} + 8 \alpha^{12} \lambda_0^5 \lambda_2
      (3 \lambda_1^2+ 22 \lambda_0  \lambda_2) + 64 \alpha^2 \lambda_2^6
      (9 \lambda_1^2 + 68 \lambda_0 \lambda_2) + 32 \alpha^4 \lambda_2^4
      (3 \lambda_1^4 + 44 \lambda_0 \lambda_1^2 \lambda_2\\
\qquad{} +
       200 \lambda_0^2 \lambda_2^2)+ 4 \alpha^{10} \lambda_0^3 \lambda_2
      (2 \lambda_1^4 + 47 \lambda_0 \lambda_1^2 \lambda_2
  + 220 \lambda_0^2 \lambda_2^2) + 8 \alpha^8 \lambda_0^2
\lambda_2^2
      (7 \lambda_1^4 + 88 \lambda_0 \lambda_1^2 \lambda_2\\
\qquad{} +
       340 \lambda_0^2 \lambda_2^2) + 4 \alpha^6 \lambda_2^2
      (3 \lambda_1^6
      + 32 \lambda_0 \lambda_1^4 \lambda_2 + 348 \lambda_0^2
\lambda_1^2
        \lambda_2^2 + 1328 \lambda_0^3 \lambda_2^3))
 \lambda_3 s_1 s_3^2\\
 \qquad{} +
   4 \alpha \lambda_1 \lambda_2 (\alpha^4 \lambda_0^2 + 4 \lambda_2^2 +
     \alpha^2 (\lambda_1^2 + 4 \lambda_0 \lambda_2)) (128 \lambda_2^7
     - 64 \alpha^2 \lambda_2^5 (\lambda_1^2 - 6 \lambda_0
\lambda_2)\\
\qquad{} +
     \alpha^{12} \lambda_0^5 (2 \lambda_0 \lambda_2-\lambda_1^2)
  -
     8 \alpha^4 \lambda_2^3 (\lambda_1^4 + 18 \lambda_0 \lambda_1^2 \lambda_2 -
       60 \lambda_0^2 \lambda_2^2) + \alpha^{10} \lambda_0^3
      ( 24
\lambda_0^2 \lambda_2^2 -\lambda_1^4\\
\qquad{} - 12 \lambda_0 \lambda_1^2 \lambda_2) +
     2 \alpha^8 \lambda_0^2 \lambda_2 (60 \lambda_0^2 \lambda_2^2-3 \lambda_1^4
     - 28 \lambda_0 \lambda_1^2 \lambda_2) - 2 \alpha^6 \lambda_2
      (\lambda_1^2 - 2 \lambda_0 \lambda_2)\\
\qquad{}\times (\lambda_1^4 + 8 \lambda_0 \lambda_1^2
        \lambda_2 + 80 \lambda_0^2 \lambda_2^2)) \lambda_3 s_3^3 -
   \alpha^3 \lambda_1 (\alpha^2 \lambda_0 + 2 \lambda_2)^2
    (\alpha^4 \lambda_0^2 + 4 \lambda_2^2 \\
   \qquad{}+
     \alpha^2 (\lambda_1^2
 + 4 \lambda_0 \lambda_2)) (\alpha^{10} \lambda_0^5 +
     64 \lambda_2^5 + \alpha^8 \lambda_0^3 (\lambda_1^2 + 12 \lambda_0\lambda_2) +
     8 \alpha^2 \lambda_2^3 (\lambda_1^2 + 18 \lambda_0 \lambda_2)\\
 \qquad{} +
     2 \alpha^6 \lambda_0^2 \lambda_2 (3 \lambda_1^2  + 28 \lambda_0 \lambda_2) +
     2 \alpha^4 \lambda_2 (\lambda_1^4 + 6 \lambda_0 \lambda_1^2 \lambda_2 +
       64 \lambda_0^2 \lambda_2^2)) \lambda_3^2 s_1^2 s_3^3 \\
 \qquad{}- 2 \alpha^4 \lambda_1^2 (\alpha^2 \lambda_0 + 2 \lambda_2)
    (\alpha^4 \lambda_0^2 + 4 \lambda_2^2+
     \alpha^2 (\lambda_1^2 + 4 \lambda_0 \lambda_2)) (\alpha^{10} \lambda_0^5+
     96 \lambda_2^5\\
 \qquad{}    + 8 \alpha^6 \lambda_0^2 \lambda_2
      (\lambda_1^2 + 9 \lambda_0 \lambda_2)
 + 16 \alpha^2 \lambda_2^3
      (\lambda_1^2 + 13 \lambda_0 \lambda_2) + \alpha^8 \lambda_0^3
      (\lambda_1^2 + 14 \lambda_0 \lambda_2) \\
 \qquad{}+ 2 \alpha^4 \lambda_2
      (\lambda_1^4 + 10 \lambda_0 \lambda_1^2 \lambda_2 + 88 \lambda_0^2 \lambda_2^2))
 \lambda_3^2 s_1 s_3^4+ \alpha^3 \lambda_1 (\alpha^4
\lambda_0^2 + 4 \lambda_2^2 +
     \alpha^2 (\lambda_1^2 + 4 \lambda_0 \lambda_2))\\
\qquad{}\times (128 \lambda_2^7 -
     96 \alpha^2 \lambda_2^5 (\lambda_1^2 - 4 \lambda_0 \lambda_2) +
     \alpha^{12} \lambda_0^5 (2 \lambda_0 \lambda_2-\lambda_1^2)
     + \alpha^{10} \lambda_0^3 (24 \lambda_0^2 \lambda_2^2 -\lambda_1^4 \\
\qquad{}- 14 \lambda_0 \lambda_1^2 \lambda_2 ) + 2 \alpha^8
\lambda_2
      (60 \lambda_0^2 \lambda_2^2-5 \lambda_1^4 - 36 \lambda_0 \lambda_1^2 \lambda_2 ) + 8 \alpha^4 \lambda_2^3
      ( 60 \lambda_0^2 \lambda_2^2 -3 \lambda_1^4\\
\qquad{}      - 26 \lambda_0 \lambda_1^2 \lambda_2)
 - 2 \alpha^6 \lambda_2
      (\lambda_1^6 + 14 \lambda_0 \lambda_1^4 \lambda_2 + 88 \lambda_0^2 \lambda_1^2
        \lambda_2^2 - 160 \lambda_0^3 \lambda_2^3)) \lambda_3^2 s_3^5\\
 \qquad{} -
   \alpha^5 \lambda_1 (\alpha^2 \lambda_0 + 2 \lambda_2)^4 \alpha^4 \lambda_0^2 + 4 \lambda_2^2 + \alpha^2
(\lambda_1^2 +
        4 \lambda_0 \lambda_2))^2
 \times ( \lambda_3^3 s_1^2 s_3^5 -
   2 \alpha^6 \lambda_1^2 (\alpha^2 \lambda_0 + 2 \lambda_2)^3\!\!\\
 \qquad{}\times
    (\alpha^4 \lambda_0^2 + 4 \lambda_2^2 + \alpha^2 (\lambda_1^2 + 4 \lambda_0 \lambda_2))^2 \lambda_3^3 s_1 s_3^6 -
   \alpha^7 \lambda_1^3 (\alpha^2 \lambda_0 + 2 \lambda_2)^2
    (\alpha^4 \lambda_0^2  + 4 \lambda_2^2 \\
   \qquad{}+ \alpha^2 (\lambda_1^2 +
        4 \lambda_0 \lambda_2))^2 \lambda_3^3 s_3^7 = 0, \\
  -2 \lambda_2 (\alpha^2 \lambda_0 + 2 \lambda_2)^3
(\alpha^{16} \lambda_0^8 +
     256 \lambda_2^8 + 16 \alpha^{12} \lambda_0^5 \lambda_2
      (\lambda_1^2 + 7 \lambda_0 \lambda_2) + \alpha^{14} \lambda_0^6
      (\lambda_1^2 + 16 \lambda_0 \lambda_2)\\
\qquad{} + 64 \alpha^2 \lambda_2^6
      (3 \lambda_1^2
 + 16 \lambda_0 \lambda_2) + 32 \alpha^4 \lambda_2^4
      (\lambda_1^4 + 16 \lambda_0 \lambda_1^2 \lambda_2 + 56 \lambda_0^2 \lambda_2^2) +
     4 \alpha^{10} \lambda_0^3 \lambda_2 (\lambda_1^4\\
\qquad{} + 25 \lambda_0 \lambda_1^2 \lambda_2 +
       112 \lambda_0^2 \lambda_2^2) + 8 \alpha^8 \lambda_0^2 \lambda_2^2
      (3 \lambda_1^4  + 40 \lambda_0 \lambda_1^2 \lambda_2 +
       140 \lambda_0^2 \lambda_2^2) + 4 \alpha^6 \lambda_2^2
      (\lambda_1^6\\
\qquad{} + 12 \lambda_0 \lambda_1^4 \lambda_2 + 140 \lambda_0^2
\lambda_1^2
        \lambda_2^2 + 448 \lambda_0^3 \lambda_2^3)) r_1 -
   2 \lambda_1 \lambda_2 (\alpha^2 \lambda_0 + 2 \lambda_2)^2 (\alpha^{16} \lambda_0^8
   + 256 \lambda_2^8 \\
\qquad{}+ 16 \alpha^{12} \lambda_0^5 \lambda_2
      (\lambda_1^2 + 7 \lambda_0 \lambda_2) + \alpha^{14} \lambda_0^6
      (\lambda_1^2 + 16 \lambda_0 \lambda_2) + 64 \alpha^2 \lambda_2^6
      (3 \lambda_1^2 + 16 \lambda_0 \lambda_2) \\
\qquad{}+ 32 \alpha^4
      \lambda_2^4    (\lambda_1^4 + 16 \lambda_0 \lambda_1^2 \lambda_2 + 56
\lambda_0^2 \lambda_2^2) +
     4 \alpha^{10} \lambda_0^3 \lambda_2 (\lambda_1^4 + 25 \lambda_0 \lambda_1^2 \lambda_2 +
       112 \lambda_0^2 \lambda_2^2)\\
\qquad{} + 8 \alpha^8 \lambda_0^2 \lambda_2^2
      (3 \lambda_1^4 + 40 \lambda_0 \lambda_1^2 \lambda_2  + 140 \lambda_0^2 \lambda_2^2) + 4 \alpha^6 \lambda_2^2
      (\lambda_1^6 + 12 \lambda_0 \lambda_1^4 \lambda_2
 + 140 \lambda_0^2 \lambda_1^2
        \lambda_2^2\\
\qquad{} + 448 \lambda_0^3 \lambda_2^3)) s_1 +
   2 \alpha \lambda_2 (\alpha^2 \lambda_0 + 2 \lambda_2) (512 \lambda_0 \lambda_2^9 +
     \alpha^{16} \lambda_0^8
  (2 \lambda_0 \lambda_2 -\lambda_1^2) -
     128 \alpha^2 \lambda_2^6 (\lambda_1^4\\
  \qquad{} + \lambda_0 \lambda_1^2 \lambda_2 -
       16 \lambda_0^2 \lambda_2^2) + \alpha^{14} \lambda_0^6
      (32 \lambda_0^2 \lambda_2^2-\lambda_1^4 - 16 \lambda_0 \lambda_1^2 \lambda_2) +
     2 \alpha^{12} \lambda_0^5 \lambda_2 (112 \lambda_0^2 \lambda_2^2 \\
\qquad{} -9 \lambda_1^4 - 52 \lambda_0 \lambda_1^2
        \lambda_2 ) - 32 \alpha^4 \lambda_2^4
      (\lambda_1^6 + 13 \lambda_0 \lambda_1^4 \lambda_2 + 16 \lambda_0^2 \lambda_1^2
        \lambda_2^2 - 112 \lambda_0^3 \lambda_2^3) \\
\qquad{}+ 4 \alpha^{10} \lambda_0^3 \lambda_2
      (224 \lambda_0^3 \lambda_2^3 - 90 \lambda_0^2 \lambda_1^2
        \lambda_2^2  -\lambda_1^6 - 28 \lambda_0 \lambda_1^4 \lambda_2
       ) + 4 \alpha^8 \lambda_0^2 \lambda_2^2
      (560 \lambda_0^3 \lambda_2^3-7 \lambda_1^6\\
\qquad{} - 84 \lambda_0 \lambda_1^4 \lambda_2 - 180 \lambda_0^2
\lambda_1^2
        \lambda_2^2 ) - 4 \alpha^6 \lambda_2^2
      (\lambda_1^8 + 14 \lambda_0 \lambda_1^6 \lambda_2 + 132 \lambda_0^2
  \times \lambda_1^4 \lambda_2^2 + 208 \lambda_0^3 \lambda_1^2
\lambda_2^3 \\
\qquad{}-
       896 \lambda_0^4 \lambda_2^4)) s_3 - \alpha (\alpha^2 \lambda_0 + 2 \lambda_2)^3
    (\alpha^{16} \lambda_0^8 + 256 \lambda_2^8 + 16 \alpha^{12} \lambda_0^5 \lambda_2
      (\lambda_1^2 + 7 \lambda_0 \lambda_2)\\
\qquad{} + \alpha^{14} \lambda_0^6
      (\lambda_1^2 + 16 \lambda_0 \lambda_2) + 64 \alpha^2 \lambda_2^6
      (3 \lambda_1^2 + 16 \lambda_0 \lambda_2) + 32 \alpha^4 \lambda_2^4
      (\lambda_1^4 + 16 \lambda_0 \lambda_1^2 \lambda_2 + 56 \lambda_0^2 \lambda_2^2)\\
\qquad{} +
     4 \alpha^{10} \lambda_0^3 \lambda_2 (\lambda_1^4 + 25 \lambda_0
 \times \lambda_1^2 \lambda_2 +
       112 \lambda_0^2 \lambda_2^2) + 8 \alpha^8 \lambda_0^2 \lambda_2^2
      (3 \lambda_1^4 + 40 \lambda_0 \lambda_1^2 \lambda_2 +
       140 \lambda_0^2 \lambda_2^2)\\
 \qquad{} + 4 \alpha^6 \lambda_2^2
      (\lambda_1^6 + 12 \lambda_0 \lambda_1^4 \lambda_2 + 140 \lambda_0^2 \lambda_1^2
        \lambda_2^2 + 448 \lambda_0^3 \lambda_2^3)) \lambda_3 s_1^2 s_3 -
   2 \alpha^4 \lambda_1 (\alpha^2 \lambda_0 + 2 \lambda_2)^2\\
 \qquad{}\times
    (\alpha^{14} \lambda_0^8 + 128 \lambda_2^6 (\lambda_1^2 + 2 \lambda_0 \lambda_2) +
     \alpha^{12} \lambda_0^6 (\lambda_1^2 + 16 \lambda_0 \lambda_2) +
     4 \alpha^8 \lambda_0^3 \lambda_2 (\lambda_1^2 + 5 \lambda_0 \lambda_2)\\
  \qquad{}\times
      (\lambda_1^2 + 20 \lambda_0 \lambda_2) + 4 \alpha^{10} \lambda_0^5 \lambda_2
      (4 \lambda_1^2 + 27 \lambda_0 \lambda_2) + 32 \alpha^2 \lambda_2^4
      (\lambda_1^4 + 13 \lambda_0 \lambda_1^2 \lambda_2 + 26 \lambda_0^2 \lambda_2^2)\\
  \qquad{} +
     4 \alpha^6 \lambda_0^2 \lambda_2^2 (7 \lambda_1^4 + 78 \lambda_0 \lambda_1^2
        \lambda_2 + 220 \lambda_0^2 \lambda_2^2) + 4 \alpha^4 \lambda_2^2
      (\lambda_1^6 + 14 \lambda_0 \lambda_1^4 \lambda_2
      + 128 \lambda_0^2 \lambda_1^2
        \lambda_2^2 \\
 \qquad{}+ 288 \lambda_0^3 \lambda_2^3)) \lambda_3 s_1 s_3^2 +
   \alpha (512 \alpha^2 \lambda_2^9
      (\lambda_1^2 - 16 \lambda_0 \lambda_2)
 - 2048 \lambda_2^11 + \alpha^{20} \lambda_0^9
      (2 \lambda_0 \lambda_2-\lambda_1^2) \\
\qquad{}+ \alpha^{18} \lambda_0^7
      (2 \lambda_0 \lambda_2-\lambda_1^2) (\lambda_1^2 + 16 \lambda_0 \lambda_2) +
     256 \alpha^4 \lambda_2^7 (\lambda_1^4 + 9 \lambda_0 \lambda_1^2 \lambda_2 -
       54 \lambda_0^2 \lambda_2^2) \\
 \qquad {}+ 8 \alpha^{16} \lambda_0^6 \lambda_2
      (27 \lambda_0^2 \lambda_2^2-2 \lambda_1^4 - 9 \lambda_0 \lambda_1^2 \lambda_2
      ) +
     32 \alpha^6 \lambda_2^5 (\lambda_1^6 + 14 \lambda_0 \lambda_1^4 \lambda_2 +
       128 \lambda_0^2 \lambda_1^2 \lambda_2^2\\
\qquad{} - 384 \lambda_0^3 \lambda_2^3) +
     16 \alpha^{12} \lambda_0^3 \lambda_2^2
 \times (84 \lambda_0^3 \lambda_2^3 -2 \lambda_1^6 - 15
\lambda_0 \lambda_1^4
        \lambda_2 + 16 \lambda_0^2 \lambda_1^2 \lambda_2^2 ) - 32 \alpha^8 \lambda_0 \lambda_2^4
      (\lambda_1^6 \\
\qquad{}- 4 \lambda_0 \lambda_1^4 \lambda_2 - 116 \lambda_0^2
\lambda_1^2
        \lambda_2^2 + 168 \lambda_0^3 \lambda_2^3)
        + 4 \alpha^{14} \lambda_0^4 \lambda_2
  (192 \lambda_0^3 \lambda_2^3-\lambda_1^6 - 23 \lambda_0
\lambda_1^4 \lambda_2\\
\qquad{}      - 32 \lambda_0^2 \lambda_1^2
        \lambda_2^2 ) + 4 \alpha^{10} \lambda_0 \lambda_1^2
      \lambda_2^2 (432 \lambda_0^3 \lambda_2^3 -\lambda_1^6 - 18 \lambda_0 \lambda_1^4 \lambda_2 -
       60 \lambda_0^2 \lambda_1^2 \lambda_2^2)) \lambda_3
    s_3^3 \\
\qquad{} + 2 \alpha^3 \lambda_2 (\alpha^2 \lambda_0 + 2
\lambda_2)^4
    (\alpha^4 \lambda_0^2 + 4 \lambda_2^2 +
     \alpha^2 (\lambda_1^2 + 4 \lambda_0 \lambda_2)) (\alpha^6 \lambda_0^3 +
     12 \alpha^2 \lambda_0 \lambda_2^2 + 8 \lambda_2^3\\
\qquad{} +
     \alpha^4 \lambda_0 (6 \lambda_0 \lambda_2 -\lambda_1^2)) \lambda_3^2 s_1^2 s_3^3
      + 8 \alpha^4 \lambda_1 \lambda_2 (\alpha^2 \lambda_0 + 2
\lambda_2)
    (\alpha^4 \lambda_0^2 + 4 \lambda_2^2 +
     \alpha^2 (\lambda_1^2 + 4 \lambda_0 \lambda_2)) \\
\qquad{}\times (\alpha^{10} \lambda_0^5 +
     10 \alpha^8 \lambda_0^4 \lambda_2 + 32 \lambda_2^5 +
     8 \alpha^2 \lambda_2^3 (\lambda_1^2 + 10 \lambda_0 \lambda_2)
 + 2 \alpha^6 \lambda_0^2 \lambda_2 (\lambda_1^2 + 20 \lambda_0
\lambda_2)\\
\qquad{} +
     \alpha^4 \lambda_2 (\lambda_1^4 + 8 \lambda_0 \lambda_1^2 \lambda_2 +
       80 \lambda_0^2 \lambda_2^2)) \lambda_3^2 s_1 s_3^4 -
   2 \alpha^3 \lambda_2 (\alpha^4 \lambda_0^2 + 4 \lambda_2^2 +
     \alpha^2 (\lambda_1^2 + 4 \lambda_0 \lambda_2))
\\
\qquad{} \times (\alpha^{14} \lambda_0^7 +
     128 \lambda_2^7 - 32 \alpha^2 \lambda_2^5 (\lambda_1^2 - 14 \lambda_0 \lambda_2) +
     2 \alpha^{12} \lambda_0^5 (7 \lambda_0 \lambda_2-\lambda_1^2) -
     24 \alpha^4 \lambda_2^3 (\lambda_1^4 \\
\qquad{}+ 4 \lambda_0 \lambda_1^2 \lambda_2 -
       28 \lambda_0^2 \lambda_2^2)
 + \alpha^{10} \lambda_0^3
      (84 \lambda_0^2 \lambda_2^2-\lambda_1^4 - 18 \lambda_0 \lambda_1^2 \lambda_2
      ) +
     2 \alpha^8 \lambda_0^2 \lambda_2 (140 \lambda_0^2 \lambda_2^2-5 \lambda_1^4\\
\qquad{} - 32 \lambda_0 \lambda_1^2
        \lambda_2 ) - 4 \alpha^6 \lambda_2
      (\lambda_1^6 + 7 \lambda_0 \lambda_1^4 \lambda_2
 + 28 \lambda_0^2 \lambda_1^2
        \lambda_2^2 - 140 \lambda_0^3 \lambda_2^3)) \lambda_3^2 s_3^5 +
   \alpha^5 (\alpha^2 \lambda_0 + 2 \lambda_2)^2\\
\qquad{}\times
    (\alpha^4 \lambda_0^2 + 4 \lambda_2^2 + \alpha^2 (\lambda_1^2 +
        4 \lambda_0 \lambda_2))^2 (\alpha^6 \lambda_0^3 +
     6 \alpha^4 \lambda_0^2 \lambda_2 + 8 \lambda_2^3 \\
\qquad{}
 + 2 \alpha^2 \lambda_2
      (\lambda_1^2 + 6 \lambda_0 \lambda_2)) \lambda_3^3 s_1^2 s_3^5 +
   2 \alpha^6 \lambda_1 (\alpha^2 \lambda_0 + 2 \lambda_2)
    (\alpha^4 \lambda_0^2 + 4 \lambda_2^2 + \alpha^2 (\lambda_1^2 +
        4 \lambda_0 \lambda_2))^2 \\
\qquad{}\times(\alpha^6 \lambda_0^3 +
     6 \alpha^4 \lambda_0^2 \lambda_2 + 8 \lambda_2^3
 + 2 \alpha^2 \lambda_2
      (\lambda_1^2 + 6 \lambda_0 \lambda_2)) \lambda_3^3 s_1 s_3^6 +
   \alpha^7 \lambda_1^2 (\alpha^4 \lambda_0^2 + 4 \lambda_2^2\\
 \qquad{} +
      \alpha^2 (\lambda_1^2 + 4 \lambda_0 \lambda_2))^2 (\alpha^6 \lambda_0^3 +
     6 \alpha^4 \lambda_0^2 \lambda_2 + 8 \lambda_2^3 + 2 \alpha^2 \lambda_2
      (\lambda_1^2
 + 6 \lambda_0 \lambda_2)) \lambda_3^3 s_3^7 = 0.
\end{gather*}

\vspace{-3mm}

\section{Appendix}
 The stationary solutions of the system (\ref{005511}):
\begin{gather*}
  \Bigg \{ \Bigg\{
  s_1 = -\frac{\sqrt{ 2 \alpha^2 \lambda_0^2 + \lambda_1^2
  - \lambda_1 \sqrt{4 \alpha^2 \lambda_0^2 + \lambda_1^2}} \
  \big(\lambda_1 + \sqrt{4 \alpha^2 \lambda_0^2 + \lambda_1^2}\big)
  \sqrt{\lambda_0 - \lambda_3 s_2^2}}
  {2 \sqrt{2} \alpha^2 \lambda_0^2 \sqrt{\lambda_3}}, \nonumber \\
 s_3 = 0, \ \ r_1 = -\frac{\sqrt{ 2 \alpha^2 \lambda_0^2 +
\lambda_1^2
  - \lambda_1 \sqrt{4 \alpha^2 \lambda_0^2 + \lambda_1^2}}
  \sqrt{\lambda_0 - \lambda_3 s_2^2}}
  {2 \sqrt{2} \alpha^2 \lambda_0 \sqrt{\lambda_3}}, \\
 r_2 = -\frac{\big(\lambda_1 - \sqrt{4 \alpha^2 \lambda_0^2 +
\lambda_1^2}) s_2}
    {2 \alpha^2 \lambda_0}, \ \  r_3 = 0  \Bigg\}, \\
 \Bigg\{
  s_1 = \frac{\sqrt{ 2 \alpha^2 \lambda_0^2 + \lambda_1^2
  - \lambda_1 \sqrt{4 \alpha^2 \lambda_0^2 + \lambda_1^2}}
  \big(\lambda_1 + \sqrt{4 \alpha^2 \lambda_0^2 + \lambda_1^2}\big)
  \sqrt{\lambda_0 - \lambda_3 s_2^2}}
  {2 \sqrt{2} \alpha^2 \lambda_0^2 \sqrt{\lambda_3}}, \nonumber \\
 s_3 = 0, \ \  r_1 = \frac{\sqrt{ 2 \alpha^2 \lambda_0^2 +
\lambda_1^2
  - \lambda_1 \sqrt{4 \alpha^2 \lambda_0^2 + \lambda_1^2}}
  \sqrt{\lambda_0 - \lambda_3 s_2^2}}
  {2 \sqrt{2} \alpha^2 \lambda_0 \sqrt{\lambda_3}}, \\
 r_2 = -\frac{(\lambda_1 - \sqrt{4 \alpha^2 \lambda_0^2 +
\lambda_1^2}) s_2}
    {2 \alpha^2 \lambda_0}, \ \ r_3 = 0 \nonumber
 \Bigg\} \Bigg \}.
\end{gather*}
Here
 $\lambda_2 = -\big(2 \alpha^2 \lambda_0^2 + \lambda_1^2
 + \lambda_1 \sqrt{4 \alpha^2 \lambda_0^2 + \lambda_1^2}\big)/(4 \lambda_0)$.
\begin{gather*}
 \Bigg \{ \Bigg\{
  s_1 = -\frac{\sqrt{ 2 \alpha^2 \lambda_0^2 + \lambda_1^2
  + \lambda_1 \sqrt{4 \alpha^2 \lambda_0^2 + \lambda_1^2}}
  \big(\lambda_1 - \sqrt{4 \alpha^2 \lambda_0^2 + \lambda_1^2}\big)
  \sqrt{\lambda_0 - \lambda_3 s_2^2}}
  {2 \sqrt{2} \alpha^2 \lambda_0^2 \sqrt{\lambda_3}}, \nonumber \\
 s_3 = 0, \ \ r_1 = -\frac{\sqrt{ 2 \alpha^2 \lambda_0^2 +
\lambda_1^2
  + \lambda_1 \sqrt{4 \alpha^2 \lambda_0^2 + \lambda_1^2}}
  \sqrt{\lambda_0 - \lambda_3 s_2^2}}
  {2 \sqrt{2} \alpha^2 \lambda_0 \sqrt{\lambda_3}}, \\
 r_2 = -\frac{\big(\lambda_1 + \sqrt{4 \alpha^2 \lambda_0^2 +
\lambda_1^2}\big) s_2}
    {2 \alpha^2 \lambda_0}, \ \ r_3 = 0 \nonumber
 \Bigg\}, \\
 \Bigg\{
  s_1 = \frac{\sqrt{ 2 \alpha^2 \lambda_0^2 + \lambda_1^2
  + \lambda_1 \sqrt{4 \alpha^2 \lambda_0^2 + \lambda_1^2}}
  \big(\lambda_1 - \sqrt{4 \alpha^2 \lambda_0^2 + \lambda_1^2}\big)
  \sqrt{\lambda_0 - \lambda_3 s_2^2}}
  {2 \sqrt{2} \alpha^2 \lambda_0^2 \sqrt{\lambda_3}}, \nonumber \\
 s_3 = 0, \ \ r_1 = \frac{\sqrt{ 2 \alpha^2 \lambda_0^2 +
\lambda_1^2
  + \lambda_1 \sqrt{4 \alpha^2 \lambda_0^2 + \lambda_1^2}}
  \sqrt{\lambda_0 - \lambda_3 s_2^2}}
  {2 \sqrt{2} \alpha^2 \lambda_0 \sqrt{\lambda_3}}, \\
 r_2 = -\frac{(\lambda_1 + \sqrt{4 \alpha^2 \lambda_0^2 +
\lambda_1^2}) s_2}
    {2 \alpha^2 \lambda_0}, \ \ r_3 = 0 \Bigg \} \Bigg \}.
\end{gather*}
Here
 $ \lambda_2 = -\big(2 \alpha^2 \lambda_0^2 + \lambda_1^2
 - \lambda_1 \sqrt{4 \alpha^2 \lambda_0^2 + \lambda_1^2}\big)/(4 \lambda_0)$.

\LastPageEnding

\end{document}